\documentclass[12pt, english]{article}
\pdfoutput=1
\usepackage[margin=30mm]{geometry}
\usepackage{relsize}
\usepackage{amsmath}
\usepackage{amsfonts}
\usepackage{amssymb}
\usepackage{graphicx}
\usepackage{verbatim}
\immediate\write18{texcount -tex -sum  \jobname.tex > \jobname.wordcount.tex}

\usepackage{sectsty}
\sectionfont{\fontsize{14}{13}\selectfont}

\usepackage{amssymb}
\usepackage{tikz}
\usepackage{graphicx}
\usepackage{amsmath}
\usepackage{amsthm}
\usepackage{amssymb}
\usepackage{amsfonts, mathtools, amsbsy, dsfont}
\usepackage{xcolor}

\usepackage[utf8]{inputenx}
\usepackage{indentfirst}
\usepackage{microtype}

\usepackage[nottoc, notlof, notlot]{tocbibind}
\usepackage{hyperref}

\newcommand*\dif{\mathop{}\!\mathrm{d}}
\newcommand{\floor}[1]{\left\lfloor #1 \right\rfloor}

\newtheorem{theorem}{Theorem}[section] 
\newtheorem{corollary} {Corollary}[section] 
\newtheorem{proposition}{Proposition}[section] 

\theoremstyle{definition} 
\newtheorem{remark}{Remark}[section] 

\numberwithin{equation}{section} 

\providecommand{\keywords}[1]
{
  \small	
  \textbf{\textit{Keywords: }} #1
}
\providecommand{\MSC}[1]
{
  \small	
  \textit{2010 MSC: } #1   
}

\title{On the Exact Distributions of the Maximum of the Asymmetric Telegraph Process}
\author{Fabrizio Cinque$^1$, Enzo Orsingher$^2$\\
        \small Department of Statistical Sciences, Sapienza University of Rome, Italy \\
        \small $^1$ cinque.fabrizio@gmail.com, \small $^2$ enzo.orsingher@uniroma1.it
}
\date{} 

\begin{document}

\maketitle

\begin{abstract}
In this paper we present the distribution of the maximum of the asymmetric telegraph process in an arbitrary time interval $[0,t]$ under the conditions that the initial velocity $V(0)$ is either $c_1$ or $-c_2$ and the number of changes of direction is odd or even. 
For the case $V(0) = -c_2$ the singular component of the distribution of the maximum displays an unexpected cyclic behavior and depends only on $c_1$ and $c_2$, but not on the current time $t$.
We obtain also the unconditional distribution of the maximum for either $V(0) = c_1$ or $V(0) = -c_2$ and its expression has the form of series of Bessel functions.
We also show that all the conditional distributions emerging in this analysis are governed by generalized Euler-Poisson-Darboux equations.
We recover all the distributions of the maximum of the symmetric telegraph process as particular cases of the present paper. 
\\We underline that it rarely happens to obtain explicitly the distribution of the maximum of a process. For this reason the results on the range of oscillations of a natural process like the telegraph model make it useful for many applications.
\end{abstract} \hspace{10pt}


\MSC{Primary 60K99}

\section{Introduction}
We consider a particle starting from the origin and moving rightward with velocity $c_1$ and leftward with velocity $-c_2$ initially taken with equal probability. The change of velocity is paced by a homogeneous Poisson process $\{N(t)\}_{t\ge0}$, with rate $\lambda>0$. Some typical sample paths of the asymmetric telegraph process are depicted in Figure \ref{traiettoria1} and Figure \ref{traiettoria2}.

The different velocities of motion introduce asymmetry in the sample paths and in all the related distribution functions and this makes the evaluation of the distribution of the maximum much more difficult than in the symmetric case, see \cite{CO2020}.
\begin{figure}
\begin{minipage}{0.5\textwidth}
\centering
\begin{tikzpicture}[scale = 0.65]
\draw [lightgray] (7.3,2.6) -- (0,2.6) node[left, black, scale = 1]{$\beta$}; \draw[dashed, lightgray] (7.3, 2.6) -- (8.4,2.6);
\draw (0,0) -- (1,1.3)--(1.6,0.22)--(2.8, 1.78)--(3.2, 1.06)--(4.1, 2.23)--(5, 0.61)--(5.4, 1.13)--(6.8, -1.39) -- (7,-1.13) ;
\draw[dashed, gray](1,1.3)-- (1,0) node[below, black, scale = 0.8]{$t_1$};
\draw[dashed, gray] (1.6,0.22)-- (1.6,0) node[below, black, scale = 0.8]{$t_2$};
\draw[dashed, gray] (2.8, 1.78)-- (2.8,0) node[below, black, scale = 0.8]{$t_3$};
\draw[dashed, gray] (3.2, 1.06 )-- (3.2,0) node[below, black, scale = 0.8]{$t_4$};
\draw[dashed, gray] (4.1, 2.23)-- (4.1,0) node[below, black, scale = 0.8]{$t_5$};
\draw[dashed, gray] (5, 0.61)-- (5,0) node[below, black, scale = 0.8]{$t_6$};
\draw[dashed, gray] (5.4, 1.13)-- (5.4,0) node[below, black, scale = 0.8]{$t_7$};
\draw[dashed, gray] (6.8, -1.39)-- (6.8,0) node[above, black, scale = 0.8]{$t_8\ \ $};
\draw[dashed, gray] (7,-1.13)-- (7,0) node[above, black, scale = 1]{$\ t$};
\draw[->, thick] (-1,0) -- (8,0) node[below, scale = 1]{$\pmb{s}$};
\draw[->, thick] (0,-1.7) -- (0,4) node[left, scale = 1]{ $\pmb{\mathcal{T}(s)}$};
\filldraw[black] (0,0) circle (3.5pt) node[below left, black, scale =0.9]{\textbf{O}};
\end{tikzpicture}
\caption{motion with $V(0) = c_1$ and $\newline c_1<c_2$}\label{traiettoria1}
\end{minipage}
\begin{minipage}{0.5\textwidth}
\centering
\begin{tikzpicture}[scale = 0.65]
\draw [lightgray] (7.3,2.6) -- (0,2.6) node[left, black, scale = 1]{$\beta$}; \draw[dashed, lightgray] (7.3, 2.6) -- (8.4,2.6);
\draw (0,0) -- (1.1, -1.32) --(2,0.12)--(2.8,-0.84 )--(4, 1.08) --(4.4, 0.6)--(5.2, 1.88)--(5.7, 1.28)--(7, 3.36);
\draw[dashed, gray] (1.1, -1.32)-- (1.1,0) node[above, black, scale = 0.8]{$t_1$};
\draw[dashed, gray] (2,0.12)-- (2,0) node[below, black, scale = 0.8]{$t_2$};
\draw[dashed, gray] (2.8,-0.84 )-- (2.8,0) node[above, black, scale = 0.8]{$t_3$};
\draw[dashed, gray] (4, 1.08)-- (4,0) node[below, black, scale = 0.8]{$t_4\ $};
\draw[dashed, gray] (4.4, 0.6)-- (4.4,0) node[below, black, scale = 0.8]{$\ t_5$};
\draw[dashed, gray] (5.2, 1.88)-- (5.2,0) node[below, black, scale = 0.8]{$t_6$};
\draw[dashed, gray] (5.7, 1.28)-- (5.7,0) node[below, black, scale = 0.8]{$t_7$};
\draw[dashed, gray] (7, 3.36)-- (7,0) node[below, black, scale = 1]{$t$};
\draw[->, thick] (-1,0) -- (8,0) node[below, scale = 1]{$\pmb{s}$};
\draw[->, thick] (0,-1.7) -- (0,4) node[left, scale = 1]{ $\pmb{\mathcal{T}(s)}$};
\filldraw[black] (0,0) circle (3.5pt) node[below left, black, scale =0.9]{\textbf{O}};
\end{tikzpicture}
\caption{motion with $V(0) =-c_2$ and \newline$c_1>c_2$}\label{traiettoria2}
\end{minipage}
\end{figure}

As far as we know, the asymmetric telegraph process was first dealt with in a paper by \cite{C1975} where the drift, due also to different rates of reversal, is eliminated by means of the relativistic Lorentz transformations. This investigation was further continued in \cite{O1987} and the explicit distribution of the position of the asymmetric telegraph process was obtained in two different ways in the paper by \cite{BNO2001}. The derivation of the distribution of the maximum of the telegraph process and the related first passage time was carried out, in the symmetric case, in \cite{O1990}, \cite{F1992} and \cite{FK1994} and recently by a different approach and under all possible conditions on the initial velocity and on the number of changes of direction by \cite{CO2020}. In a more general case, the first passage time was presented by \cite{SZ2004} and \cite{Z2004}. Here we make a further step forward because we assume that the rightward and leftward velocities differ. This introduces a substantial generalization because the process has a drift which is of Galilean level, that is a Galilean transformation is sufficient to symmetrize the motion. We underline that the distribution of the maximum in the symmetric case displays a sort of reflection principle which is not valid in the asymmetric case treated here.
\\
We underline the fact that for few processes the explicit distribution of the maximum is known especially if a drift is assumed.
\\

The derivation of the state probability of the asymmetric telegraph process is presented in the paper by \cite{BNO2001} and the main information is that the asymmetric telegraph process $\{\mathcal{T}(t)\}_{t\ge0}$ has distribution of the absolutely continuous component, $p(x,t),\ -c_2t\le x\le c_1t$, which satisfies the hyperbolic equation
\begin{equation}\label{IntroEquazioneTelegrafoAsimmetrico}
\frac{\partial^2 p}{\partial t^2} -c_1c_2\frac{\partial^2 p}{\partial x^2} + (c_1-c_2)\frac{\partial^2 p}{\partial x\partial t} = -2\lambda \frac{\partial p}{\partial t} + \lambda(c_2-c_1) \frac{\partial p}{\partial x} 
\end{equation}
and has the explicit form
\begin{equation}\label{distribuzioneTelegrafo}
p(x,t) = \frac{e^{-\lambda t}}{c_1+c_2}\Biggl[\lambda I_0\Bigl( \frac{2\lambda}{c_1+c_2}\sqrt{(c_1t-x)(c_2t+x)}\Bigr)+
\end{equation}
$$ + \frac{\partial}{\partial t}I_0\Bigl( \frac{2\lambda}{c_1+c_2}\sqrt{(c_1t-x)(c_2t+x)}\Bigr) + \frac{c_1-c_2}{2} \frac{\partial}{\partial x}I_0\Bigl( \frac{2\lambda}{c_1+c_2}\sqrt{(c_1t-x)(c_2t+x)}\Bigr) \Biggr]$$
where $I_0(x) = \sum_{j=0}^\infty \Bigl(\frac{x}{2} \Bigr)^{2j}\frac{1}{j!^2}$, for $x \in \mathbb{R}$, is the Bessel function of the first type with imaginary argument and order $0$.
\\Clearly at both $x = -c_2t$ and $x = c_1t$ a mass of probability $\frac{e^{-\lambda t}}{2}$ is concentrated.
\\From (\ref{distribuzioneTelegrafo}) we can extract the following conditional distributions, for $k \in \mathbb{N}_0$
\begin{equation}\label{distribuzioneTelegrafoCondizionatoDispari}
P\{\mathcal{T}(t)\in \dif x\ |\ N(t) = 2k+1\}= \frac{(2k+1)!}{k!^2}\frac{(c_1t -x)^k(c_2t+x)^k}{\bigl[ (c_1+c_2)t\bigr]^{2k+1}}\dif x 
\end{equation}
and for $k\in \mathbb{N}$
\begin{equation}\label{distribuzioneTelegrafoCondizionatoPari}
P\{\mathcal{T}(t)\in \dif x\ |\ N(t) = 2k\}= \frac{(2k-1)!}{(k-1)!^2}\frac{(c_1t -x)^{k-1}(c_2t+x)^{k-1}}{\bigl[ (c_1+c_2)t\bigr]^{2k-1}}\dif x =
\end{equation}
$$ = P\{\mathcal{T}(t)\in \dif x\ |\ N(t) = 2k-1\}$$
for $-c_2t\le x\le c_1t$.
\\Further conditional relationships can be extracted when the initial velocity is fixed. For $N(t) = 2k+1$ we have $k+1$ upwards steps and $k+1$ downwards steps for both initial velocities. The exchangeability of the displacements permit us to conclude that the following holds, for all natural $k\ge0$
\begin{equation}\label{distribuzioneTelegrafoCondizionatoDispari+-}
P\{\mathcal{T}(t)\in \dif x\ |\ V(0) = c_1,\ N(t) = 2k+1\}=
\end{equation}
$$ =P\{\mathcal{T}(t)\in \dif x\ |\ V(0) = -c_2,\ N(t) = 2k+1\} = P\{\mathcal{T}(t)\in \dif x\ |\ N(t) = 2k+1\}$$
and, from (\ref{distribuzioneTelegrafoCondizionatoDispari+-}), by means of recurrence arguments, we obtain that for $k\in \mathbb{N}$ and $-c_2t\le x \le c_1t$
\begin{equation}\label{distribuzioneTelegrafoCondizionatoPari+}
P\{\mathcal{T}(t)\in \dif x\ |\ V(0) = c_1,\ N(t) = 2k\}= \frac{(2k)!}{k!(k-1)!}\frac{(c_1t -x)^{k-1}(c_2t+x)^{k}}{\bigl[ (c_1+c_2)t\bigr]^{2k}}\dif x 
\end{equation}
\begin{equation}\label{distribuzioneTelegrafoCondizionatoPari-}
P\{\mathcal{T}(t)\in \dif x\ |\ V(0) = -c_2,\ N(t) = 2k\}= \frac{(2k)!}{k!(k-1)!}\frac{(c_1t -x)^{k}(c_2t+x)^{k-1}}{\bigl[ (c_1+c_2)t\bigr]^{2k}}\dif x 
\end{equation}

A Galilean transformation 
\begin{equation}
x' = x+ \frac{c_2-c_1}{2}t\ ,\ \ \ \ \ t' = t
\end{equation}
reduces (\ref{IntroEquazioneTelegrafoAsimmetrico}) to the standard telegraph equation and (\ref{distribuzioneTelegrafo}), (\ref{distribuzioneTelegrafoCondizionatoDispari}) and (\ref{distribuzioneTelegrafoCondizionatoPari}) to the corresponding distributions of the symmetric telegraph process, see \cite{DiGOS2004}.
\\

The main results of this paper are the conditional and unconditional distributions of the maximum of the asymmetric telegraph process.
\\We show that, for $n\in \mathbb{N}$ and $\beta \in [0, c_1t]$
\begin{equation}\label{IntroMassimoRip+n}
P\{\max_{0\le s\le t}\mathcal{T}(s)< \beta\ |\ N(t) = n,\ V(0) = c_1\} = 
\end{equation}
$$ = \frac{1}{\bigl[(c_1+c_2)t \bigr]^{n}} \sum_{j=0}^{\floor{\frac{n-1}{2}}} \binom{n}{j} \Bigl[ (c_1t-\beta)^{j}(c_2t+\beta)^{n-j}  - \Bigl(\frac{c_2}{c_1}\Bigr)^{n-2j}(c_1t-\beta)^{n-j}(c_2t+\beta)^{j} \Bigr]$$
which for $c_1 = c_2 = c$ reduces to, for $n = 2k+1$
\begin{equation}\label{}
P\{\max_{0\le s\le t}\mathcal{T}(s)< \beta\ |\ V(0) = c,\ N(t) = 2k+1\} = \frac{\beta}{ct}\sum_{j=0}^{k}\binom{2j}{j}\Bigl(\frac{\sqrt{c^2t^2-\beta^2}}{2ct} \Bigr)^{2j} 
\end{equation}
with integer $k\ge0$, as shown in Corollary $5.1$ of \cite{CO2020}.

By means of (\ref{IntroMassimoRip+n}) we arrive at the distribution of the maximum of the asymmetric telegraph process with positive initial velocity, for $0\le\beta\le c_1t$
\begin{equation}\label{IntroMassimoRip+}
P\{\max_{0\le s\le t}\mathcal{T}(s)< \beta\ |\ V(0) = c_1\} = 
\end{equation}
$$ =e^{-\lambda t}\sum_{r=1}^\infty I_r\Bigl(\frac{2\lambda}{c_1+c_2} \sqrt{(c_1t-\beta)(c_2t+\beta)}\Bigr)\Biggl[ \Biggl(\sqrt{\frac{c_2t+\beta}{c_1t-\beta}}\Biggr)^r-\Biggl(\frac{c_2}{c_1}\sqrt{\frac{c_1t-\beta}{c_2t+\beta}}\Biggr)^r\ \Biggr]$$
where, for $x,r \in \mathbb{R}$
$$ I_r(x) = \sum_{j=0}^\infty \Bigl(\frac{x}{2} \Bigr)^{2j+r}\frac{1}{j!\Gamma(j+r+1)}$$
is the Bessel function of the first type with imaginary argument and order $r$.
\\We observe that at $\beta = c_1t$, formula (\ref{IntroMassimoRip+}) yields
$$P\{\max_{0\le s\le t}\mathcal{T}(s)=c_1t\ |\ V(0) = c_1\} = 1-e^{-\lambda t} $$
because the particle runs to $x = c_1t$ if no Poisson event occurs, therefore with probability $e^{-\lambda t} $.
\\

If $V(0) = -c_2$, there is a substantial difference with respect to the previous case because with a positive probability the maximum of $\mathcal{T}$ will be zero for the time interval $[0,t]$. In this case we are able to prove that, for $k\in \mathbb{N}$
\begin{equation}\label{IntroSingolarita0prima}
P\{\max_{0\le s\le t}\mathcal{T}(s)=0\ |\ V(0) = -c_2,\ N(t) = 2k-1\} =P\{\max_{0\le s\le t}\mathcal{T}(s)=0\ |\ V(0) = -c_2,\ N(t) = 2k\} =
\end{equation}
\begin{equation}\label{IntroSingolarita0}
= \sum_{j=0}^k \Biggl[\binom{2k}{j}-\binom{2k}{j-1} \Biggr] \Bigl(\frac{c_1}{c_1+c_2}\Bigr)^{j}\Bigl(\frac{c_2}{c_1+c_2}\Bigr)^{2k-j}
\end{equation}
For $c_1=c_2 =c$ we obtain that
\begin{equation}
P\{\max_{0\le s\le t}\mathcal{T}(s)=0\ |\ V(0) = -c,\ N(t) = 2k-1\} =
\end{equation}
\begin{equation}\label{IntroSingolarita0Simmetrico}
=P\{\max_{0\le s\le t}\mathcal{T}(s)=0\ |\ V(0) = -c,\ N(t) = 2k\} =\binom{2k}{k}\frac{1}{2^{2k}}
\end{equation}
which is formula $(4.4)$ of \cite{CO2020}. The probability (\ref{IntroSingolarita0prima}) does not depend on $t$, but only on $k$, as well as (\ref{IntroSingolarita0Simmetrico}). Furthermore we can represent the probabilities (\ref{IntroSingolarita0prima}) in an alternative way as, for integer $k\ge0$
\begin{equation}\label{IntroSingolarita0A}
P\{\max_{0\le s\le t}\mathcal{T}(s)=0\ |\ V(0) = -c_2,\ N(t) = 2k+1\} = \Bigl(\frac{c_2}{c_1+c_2} \Bigr)^{k+1}\sum_{j=0}^k A_j^{(k)} \Bigl(\frac{c_1}{c_1+c_2} \Bigr)^{j}
\end{equation}
where the numbers $A_j^{(k)}$ are related by the recurrence relationship
\begin{equation}
A_0^{(0)}=1,\ \ \ A_k^{(k)} = A_{k-1}^{(k)},\ \ \ A_j^{(k)} = \sum_{i=0}^j A_i^{(k-1)},\ \ \ \ k>j\ge0
\end{equation}
The numbers $A_j^{(k)}$ appearing in (\ref{IntroSingolarita0A}) form a triangular matrix which will be given below.
\\We obtain also the distribution of the maximum for an initial negative velocity. For integer $k\ge0,\ 0\le \beta\le c_1t$
\begin{equation}\label{IntroMassimoRip-p}
P\{\max_{0\le s\le t}\mathcal{T}(s)\le \beta\ |\ V(0) = -c_2,\ N(t) = 2k\} = 
\end{equation}
$$ =  \sum_{j=0}^{k} \binom{2k}{j} \frac{ (c_1t-\beta)^{j}(c_2t+\beta)^{2k-j} }{\bigl[(c_1+c_2)t \bigr]^{2k}} -\sum_{j=0}^{k-1}\binom{2k}{j} \frac{\bigl(\frac{c_2}{c_1}\bigr)^{2k-2j-1}(c_1t-\beta)^{2k-j}(c_2t+\beta)^{j}}{\bigl[(c_1+c_2)t \bigr]^{2k}} $$
In Theorem $3.2$ we obtain the distribution of the maximum under the condition that $\{V(0) = -c_2,\ N(t) = 2k+1\}$.
\\All these results permit us to obtain the unconditional distribution of the maximum of the asymmetric telegraph process with negative initial velocity in terms of series of Bessel functions. For $0\le\beta\le c_1t$
\begin{equation}\label{IntroMassimoRip-}
P\{\max_{0\le s\le t}\mathcal{T}(s)< \beta\ |\ V(0) =-c_2\} = 
\end{equation}
$$ =e^{-\lambda t}\Biggl[\sum_{r=0}^\infty I_r\Bigl(\frac{2\lambda}{c_1+c_2} \sqrt{(c_1t-\beta)(c_2t+\beta)}\Bigr)\Biggl(\sqrt{\frac{c_2t+\beta}{c_1t-\beta}}\Biggr)^r+$$
$$ -\ \frac{c_1}{c_2}\sum_{r=2}^\infty I_r\Bigl(\frac{2\lambda}{c_1+c_2} \sqrt{(c_1t-\beta)(c_2t+\beta)}\Bigr)\Biggl(\frac{c_2}{c_1}\sqrt{\frac{c_1t-\beta}{c_2t+\beta}}\Biggr)^r \ \Biggr]$$
For $\beta = 0, \ c_1=c_2=c$, we extract from (\ref{IntroMassimoRip-}) that
\begin{equation}
P\{\max_{0\le s\le t}\mathcal{T}(s)=0\ |\ V(0) =-c\} = e^{-\lambda t}\Bigl[ I_0(\lambda t)+I_1(\lambda t) \Bigr]
\end{equation}
as shown in \cite{CO2020}.
\\

We note that the conditional distributions (\ref{distribuzioneTelegrafoCondizionatoDispari+-}), (\ref{distribuzioneTelegrafoCondizionatoPari+}), (\ref{distribuzioneTelegrafoCondizionatoPari-}), as well as the terms of the distributions of the maximum, contain the functions, for integer $m,n$
\begin{equation}\label{IntroFunzionePerEPDgen}
h(x,t)= \frac{(c_1t-x)^m(c_2t+x)^n}{t^{m+n+1}}, \ \ \ \ \ -c_2t\le x \le c_1t
\end{equation}
which are strictly related to the famous generalized Euler-Poisson-Darboux equation since (\ref{IntroFunzionePerEPDgen}) solves
\begin{equation}\label{IntroEPDgen}
\frac{\partial^2 h}{\partial t^2} -c_1c_2\frac{\partial^2 h}{\partial x^2} + (c_1-c_2)\frac{\partial^2 h}{\partial x\partial t} =-\frac{m+n+2}{t}\frac{\partial h}{\partial t} -\frac{1}{t}\Bigl[ (c_1-c_2)(m+n+1)-(c_1m-c_2n)  \Bigr]\frac{\partial h}{\partial x} 
\end{equation}
which for $c_1=c_2=c$ reduces to the simplified form
\begin{equation}\label{IntroEPD}
\frac{\partial^2 h}{\partial t^2} -c^2\frac{\partial^2 h}{\partial x^2}= -\frac{m+n+2}{t}\frac{\partial h}{\partial t}+\frac{c(m-n)}{t}\frac{\partial h}{\partial x} 
\end{equation}
As a byproduct of our analysis we obtain the explicit distribution of the position of an asymmetric telegraph process where the reversals of velocity are paced by a non-homogeneous Poisson process with rate $\lambda(t)= \frac{\alpha}{t}, \ \alpha>0$, which reads, $-c_2t\le x\le c_1t$
\begin{equation}\label{}
p(x,t)=\frac{\Gamma(2\alpha)}{\Gamma(\alpha)^2}\frac{(c_1t -x)^{\alpha-1}(c_2t+x)^{\alpha-1}}{\bigl[ (c_1+c_2)t\bigr]^{2\alpha-1}}
\end{equation}
For $c_1=c_2=c$ this coincides with formula ($2.7$) of \cite{GO2016}.
\bigskip

\section{Asymmetric telegraph process with positive initial velocity}
We consider in this section the asymmetric telegraph process with positive initial velocity $V(0) = c_1$. Let us first assume $N(t) = 2k, k\in \mathbb{N}$, Poisson events.
\\Under these conditions, in order that the event $\{\max_{0\le s\le t}\mathcal{T}(s)< \beta\},\ \beta>0$, occurs, the following events must happen simultaneously
\begin{equation}\label{condizioniMassimoPari+}
\begin{cases}
c_1T_1<\beta\\
c_1T_1-c_2(T_2-T_1)+c_1(T_3-T_2)<\beta\\
\cdot\\
c_1T_1-c_2(T_2-T_1)+\dots+c_1(t-T_{2k})<\beta
\end{cases}
\end{equation}
Conditions (\ref{condizioniMassimoPari+}) can be written in a more compact form as
\begin{equation}\label{}
\bigcap_{j=0}^k \{c_1T_1-c_2(T_2-T_1)+\dots +c_1(T_{2j+1}-T_{2j})<\beta\}
\end{equation}
where $T_j$ are the random instants of occurrence of the Poisson events and coincide with the changes of direction of motion. We also assume that $T_0 = 0 \ a.s., \ T_{2k+1} = t\ a.s.$. We denote by $t_j$ the realization of the random time $T_j$.
\\
We first note that, for $k\in \mathbb{N}$
\begin{equation}\label{formulaIntegraleMassimo+p}
P\{\max_{0\le s\le t}\mathcal{T}(s)< \beta\ |\ V(0) = c_1,\ N(t) = 2k\} =
\end{equation}
$$=\sum_{j=1}^k \int_0^{\frac{\beta}{c_1}}\dif t_1 \int_{t_1}^{\frac{c_1t-\beta}{c_1+c_2}+t_1} \dif t_2 \int_{t_2}^{\frac{\beta}{c_1}+\frac{c_1+c_2}{c_1}(t_2-t_1)}\dif t_3 \int_{t_3}^{\frac{c_1t-\beta}{c_1+c_2}+t_1-t_2+t_3} \dif t_4 \dots$$
$$\dots \int_{t_{2j-2}}^{\frac{\beta}{c_1}+\frac{c_1+c_2}{c_1}(t_{2j-2}-t_{2j-3}+\dots+t_2-t_1)}\dif t_{2j-1} \int_{\frac{c_1t-\beta}{c_1+c_2}+t_1-t_2+\dots-t_{2j-2}+t_{2j-1}}^t  \frac{(2k)!}{t^{2k}}\frac{(t-t_{2j})^{2k-2j}}{(2k-2j)!}\dif t_{2j}$$
where we considered the explicit form of $P\{T_1\in \dif t_1,\dots, T_{2j}\in \dif t_{2j}|N(t) = 2k\}$.
Formula (\ref{formulaIntegraleMassimo+p}) can be written down by considering that after the obvious first condition $c_1t_1<\beta$, at each leftward step we have two possibilities. The first one is that we went so far in the left direction that we have not enough time to overcome the threshold $\beta$. The second case occurs when we moved leftward in such a way that we can reach level $\beta$, in the remaining time interval, with positive probability. For example at time $T_2 = t_2$ we can have two possibilities
\begin{equation}\label{2PossibiliCondizioniInT2}
\beta-[c_1t_1-c_2(t_2-t_1)]\le c_1(t-t_2) \ \ \ or \ \ \ \beta-[c_1t_1-c_2(t_2-t_1)]> c_1(t-t_2)
\end{equation}
In the first case we have enough time that in $(t_2,t)$ we can overcome level $\beta$ while in the second case the particle will remain below $\beta$ with probability one in the $[t_2,t]$.
\\This reasoning can be repeated at all even times $T_{2j} = t_{2j},\ j\in \mathbb{N}$. 
\\The last integral in (\ref{formulaIntegraleMassimo+p}) corresponds to the case where the particle at time $t_{2j}$ has gone so far to the left that in the interval $[t_{2j}, t]$ it will never overcome $\beta$.
\\The first condition of (\ref{2PossibiliCondizioniInT2}) gives the limits of $t_2$ while the conditions on $t_3$ is implied by
$$c_1t_1-c_2(t_2-t_1)+c_1(t_3-t_2)<\beta$$
Similar inequalities can be written for all odd times $t_{2j-1}$.
\\

In principle the distribution (\ref{formulaIntegraleMassimo+p}) can be derived by direct calculation of the integrals and for small values of $k$ this is really possible. For example, for $k=2$ we obtained that
\begin{equation}\label{esempioMassimoRip4}
P\{\max_{0\le s\le t}\mathcal{T}(s)< \beta\ |\ V(0) = c_1,\ N(t) = 4\} =  \frac{1}{\bigl[(c_1+c_2)t \bigr]^{4}}\cdot
\end{equation}
$$\cdot\Biggl[4(c_1t-\beta)(c_2t+\beta)^3-4\Bigl( \frac{c_2}{c_1}\Bigr)^2(c_1t-\beta)^3(c_2t+\beta)+(c_2t+\beta)^4- \Bigl( \frac{c_2}{c_1}\Bigr)^4(c_1t-\beta)^4 \Biggr]$$
with $\beta \in [0, c_1t]$. This and other similar calculations suggested the general expression of the distribution of the maximum which we proved by recurrence arguments as shown in the forthcoming theorem. In particular, it is necessary to keep in mind that
\begin{equation}\label{massimo+2}
P\{\max_{0\le s\le t}\mathcal{T}(s)< \beta\ |\ V(0) = c_1,\ N(t) = 2\} = \frac{(c_2t+\beta)^2-\bigl(\frac{c_2}{c_1}\bigr)^2(c_1t-\beta)^2}{(c_1+c_2)^2t^2}=
\end{equation}
$$ = \frac{\beta\bigl[\beta(c_1-c_2)+2c_1c_2t \bigr]}{(c_1+c_2)c_1^2t^2}$$
Probability (\ref{massimo+2}) is the induction basis of the recurrence method we use in the next theorem and it can be easily obtained by evaluating (\ref{formulaIntegraleMassimo+p}) for $k=1$.

\begin{theorem}\label{teoremaMassimoRip+p}
Let $\{\mathcal{T}(t)\}_{t\ge0}$ be the asymmetric telegraph process. For $k\in \mathbb{N},\ \beta\in[0, c_1t]$
\begin{equation}\label{massimoRip+p}
P\{\max_{0\le s\le t}\mathcal{T}(s)< \beta\ |\ V(0) = c_1,\ N(t) = 2k\} = 
\end{equation}
$$ = \frac{1}{\bigl[(c_1+c_2)t \bigr]^{2k}} \sum_{j=0}^{k-1} \binom{2k}{j} \Bigl[ (c_1t-\beta)^{j}(c_2t+\beta)^{2k-j}  - \Bigl(\frac{c_2}{c_1}\Bigr)^{2k-2j}(c_1t-\beta)^{2k-j}(c_2t+\beta)^{j} \Bigr]$$
\end{theorem}

We point out that the distribution (\ref{massimoRip+p}), as well as all the other distributions conditioned on $N(t)$ we are presenting below, do not depend on the rate $\lambda$ that influences the changes of direction of the motion.
\paragraph{Proof.}
The discussion above justifies that 
\begin{equation}\label{massimoRipIntegrale+p}
P\{\max_{0\le s\le t}\mathcal{T}(s)< \beta\ |\ V(0) = c_1,\ N(t) = 2k\} = 
\end{equation}
$$=\int_0^{\frac{\beta}{c_1}} \int_{t_1}^{\frac{c_1t-\beta}{c_1+c_2}+t_1} P\{\max_{0\le s\le t-t_2}\mathcal{T}(s)< \beta-(c_1+c_2)t_1+c_2t_2\ |\ V(0) = c_1,\ N(t-t_2) = 2k-2\} \cdot $$
$$ \cdot P\{T_1\in \dif t_1,\ T_2\in \dif t_2\ |\ N(t) = 2k\} +\int_0^{\frac{\beta}{c_1}} \int_{\frac{c_1t-\beta}{c_1+c_2}+t_1}^t P\{T_1\in \dif t_1,\ T_2\in \dif t_2\ |\ N(t) = 2k\}$$
The first integral in (\ref{massimoRipIntegrale+p}) is constructed by assuming that
\begin{equation}\label{}
\begin{cases}
c_1t_1<\beta\\
\beta-[c_1t_1-c_2(t_2-t_1)]\le c_1(t-t_2)
\end{cases}
\end{equation}
where the second inequality is satisfied if the moving particle has not moved leftward too much in the interval $(t_1,t_2)$ so that it has sufficient time in $(t_2,t)$ to overcome level $\beta$ with positive probability. Furthermore, we used the homogeneity of the telegraph process at the Poisson times to build the probability in the first integral.
\\If $\beta-[c_1t_1-c_2(t_2-t_1)]> c_1(t-t_2)$ the particle will remain below $\beta$ for the whole time interval $(t_2,t)$ with probability one. This explains the structure of the second integral of (\ref{massimoRipIntegrale+p}).
\\We proceed by induction by using (\ref{massimoRip+p}). We can write 
$$P\{\max_{0\le s\le t}\mathcal{T}(s)< \beta\ |\ V(0) = c_1,\ N(t) = 2k\}= \int_0^{\frac{\beta}{c_1}} \dif t_1 \int_{t_1}^{\frac{c_1t-\beta}{c_1+c_2}+t_1} \frac{(2k)!}{(2k-2)!}\frac{(t-t_2)^{2k-2}}{t^{2k}} \cdot$$
$$ \cdot\sum_{j=0}^{k-2} \frac{\binom{2k-2}{j}}{\bigl[(c_1+c_2)(t-t_2) \bigr]^{2k-2}}  \Biggl[ \bigl[c_1(t-t_2)-(\beta-(c_1+c_2)t_1+c_2t_2)\bigr]^{j}\cdot$$
$$ \cdot \bigl[c_2(t-t_2)+\beta-(c_1+c_2)t_1+c_2t_2\bigr]^{2k-2-j} - \Bigl(\frac{c_2}{c_1}\Bigr)^{2k-2-2j}\bigl[c_1(t-t_2)-(\beta-(c_1+c_2)t_1+c_2t_2)\bigr]^{2k-2-j}\cdot$$
\begin{equation}\label{passaggioIntegraleDimostrazioneMassimoRip+p}
\cdot\bigl[c_2(t-t_2)\ +\ \beta-(c_1+c_2)t_1+c_2t_2\bigr]^{j} \Biggr]\dif t_2+\Bigl(\frac{c_2t+\beta}{(c_1+c_2)t}\Bigr)^{2k}-\Bigl(\frac{c_2}{c_1}\Bigr)^{2k}\Bigl(\frac{c_1t-\beta}{(c_1+c_2)t}\Bigr)^{2k}
\end{equation}
Clearly, the last two terms of (\ref{passaggioIntegraleDimostrazioneMassimoRip+p}) represent the second integral of (\ref{massimoRipIntegrale+p}). The first integral of (\ref{passaggioIntegraleDimostrazioneMassimoRip+p}) can be substantially simplified and becomes
$$\int_0^{\frac{\beta}{c_1}} \dif t_1 \int_{t_1}^{\frac{c_1t-\beta}{c_1+c_2}+t_1} \sum_{j=0}^{k-2} \binom{2k-2}{j}\frac{(2k)!}{(2k-2)!}\frac{1}{(c_1+c_2)^{2k-2}t^{2k}}\cdot$$
$$\cdot  \Biggl[\bigl[c_1t-\beta +(c_1+c_2)(t_1-t_2) \bigr]^{j}\bigl[c_2t+\beta-(c_1+c_2)t_1 \bigr]^{2k-2-j} +$$
$$ - \Bigl(\frac{c_2}{c_1}\Bigr)^{2k-2-2j} \bigl[c_1t-\beta +(c_1+c_2)(t_1-t_2) \bigr]^{2k-2-j}\bigl[c_2t+\beta-(c_1+c_2)t_1 \bigr]^{j} \Biggr]\dif t_2=$$
$$ = \sum_{j=0}^{k-2} \binom{2k-2}{j}\frac{(2k)!}{(2k-2)!(j+1)}\frac{(c_1t-\beta)^{j+1}}{(c_1+c_2)^{2k-1}t^{2k}} \int_0^{\frac{\beta}{c_1}} \bigl[c_2t+\beta-(c_1+c_2)t_1 \bigr]^{2k-2-j} \dif t_1+$$
$$ - \sum_{j=0}^{k-2} \binom{2k-2}{j}\frac{(2k)!}{(2k-2)!(2k-1-j)}\frac{(c_1t-\beta)^{2k-1-j}}{(c_1+c_2)^{2k-1}t^{2k}} \Bigl(\frac{c_2}{c_1}\Bigr)^{2k-2-2j}\int_0^{\frac{\beta}{c_1}} \bigl[c_2t+\beta-(c_1+c_2)t_1 \bigr]^{j} \dif t_1  =$$ 
$$ = \sum_{j=0}^{k-2} \binom{2k}{j+1}\frac{(c_1t-\beta)^{j+1}}{\bigl[(c_1+c_2)t\bigr]^{2k}}\Bigl[(c_2t +\beta)^{2k-1-j}- \Bigl(\frac{c_2}{c_1}\Bigr)^{2k-1-j}(c_1t-\beta)^{2k-1-j} \Bigr]+$$
$$ -\sum_{j=0}^{k-2} \binom{2k}{j+1}\frac{(c_1t-\beta)^{2k-1-j}}{\bigl[(c_1+c_2)t\bigr]^{2k}}\Bigl(\frac{c_2}{c_1}\Bigr)^{2k-2-2j}\Bigl[(c_2t +\beta)^{j+1}- \Bigl(\frac{c_2}{c_1}\Bigr)^{j+1}(c_1t-\beta)^{j+1} \Bigr]=$$
\begin{equation}\label{ultimoPassaggioDimMassimoRip+p}
=\frac{1}{\bigl[(c_1+c_2)t \bigr]^{2k}} \sum_{j=1}^{k-1} \binom{2k}{j} \Bigl[ (c_1t-\beta)^{j}(c_2t+\beta)^{2k-j}  - \Bigl(\frac{c_2}{c_1}\Bigr)^{2k-2j}(c_1t-\beta)^{2k-j}(c_2t+\beta)^{j} \Bigr]
\end{equation}
If we add to result (\ref{ultimoPassaggioDimMassimoRip+p}) the last two terms of (\ref{passaggioIntegraleDimostrazioneMassimoRip+p}) we readily arrive at the claimed distribution (\ref{massimoRip+p}).
\\The reader can also check that for $k=2$, the result of Theorem \ref{teoremaMassimoRip+p} coincides with (\ref{esempioMassimoRip4}) obtained by evaluating directly the integral (\ref{formulaIntegraleMassimo+p}) and for $k=1$ (\ref{massimoRip+p}) coincides with (\ref{massimo+2}).
\hfill$\Box$
\\

We have now the distribution of the maximum of the asymmetric telegraph process for a rightward initial step and an odd number of changes of direction.

\begin{theorem}\label{teoremaMassimoRip+d}
Let $\{\mathcal{T}(t)\}_{t\ge0}$ be the asymmetric telegraph process. For $k\in \mathbb{N}_0,\ \beta\in [0, c_1t]$
\begin{equation}\label{massimoRip+d}
P\{\max_{0\le s\le t}\mathcal{T}(s)< \beta\ |\ V(0) = c_1,\ N(t) = 2k+1\} = 
\end{equation}
$$ = \sum_{j=0}^{k} \binom{2k+1}{j}  \frac{(c_1t-\beta)^{j}(c_2t+\beta)^{2k+1-j}  - \Bigl(\frac{c_2}{c_1}\Bigr)^{2k+1-2j}(c_1t-\beta)^{2k+1-j}(c_2t+\beta)^{j} }{\bigl[(c_1+c_2)t \bigr]^{2k+1}}$$
\end{theorem}

\paragraph{Proof.}
First of all we consider that
\begin{equation}\label{massimo+1}
P\{\max_{0\le s\le t}\mathcal{T}(s)< \beta\ |\ V(0) = c_1,\ N(t) = 1\} = P\{T_1<\frac{\beta}{c_1}\ |\ N(t) = 1\} =\frac{\beta}{c_1t}
\end{equation}
which coincides with (\ref{massimoRip+d}) for $k=1$. The distribution (\ref{massimo+1}) is the first step of our induction procedure.
\\Now we use again the recursive arguments of Theorem \ref{teoremaMassimoRip+p}.
\begin{equation}\label{massimoRipIntegrale+d}
P\{\max_{0\le s\le t}\mathcal{T}(s)< \beta\ |\ V(0) = c_1,\ N(t) = 2k+1\} = 
\end{equation}
$$=\int_0^{\frac{\beta}{c_1}} \int_{t_1}^{\frac{c_1t-\beta}{c_1+c_2}+t_1} P\{\max_{0\le s\le t-t_2}\mathcal{T}(s)< \beta-(c_1+c_2)t_1+c_2t_2\ |\ V(0) = c_1,\ N(t-t_2) = 2k-1\} \cdot $$
$$ \cdot P\{T_1\in \dif t_1,\ T_2\in \dif t_2\ |\ N(t) = 2k+1\} +\int_0^{\frac{\beta}{c_1}} \int_{\frac{c_1t-\beta}{c_1+c_2}+t_1}^t P\{T_1\in \dif t_1,\ T_2\in \dif t_2\ |\ N(t) = 2k+1\}$$
The arguments adopted in justifying (\ref{massimoRipIntegrale+p}) are the same as those needed here to write (\ref{massimoRipIntegrale+d}).
\\We must evaluate the first integral of (\ref{massimoRipIntegrale+d}) while the second one is straightforwardly determined. In view of (\ref{massimoRip+d}), that we consider our induction procedure, we can write the first integral as follows
$$\int_0^{\frac{\beta}{c_1}} \dif t_1 \int_{t_1}^{\frac{c_1t-\beta}{c_1+c_2}+t_1} \sum_{j=0}^{k-1} \binom{2k-1}{j}\frac{(2k+1)!}{(2k-1)!}\frac{1}{(c_1+c_2)^{2k-1}t^{2k+1}}\cdot$$
$$\cdot  \Biggl[\bigl[c_1t-\beta +(c_1+c_2)(t_1-t_2) \bigr]^{j}\bigl[c_2t+\beta-(c_1+c_2)t_1 \bigr]^{2k-1-j} +$$
$$ - \Bigl(\frac{c_2}{c_1}\Bigr)^{2k-1-2j} \bigl[c_1t-\beta +(c_1+c_2)(t_1-t_2) \bigr]^{2k-1-j}\bigl[c_2t+\beta-(c_1+c_2)t_1 \bigr]^{j} \Biggr]\dif t_2=$$
$$ = \sum_{j=0}^{k-1} \binom{2k-1}{j}\frac{(2k+1)!}{(2k-1)!(j+1)}\frac{(c_1t-\beta)^{j+1}}{(c_1+c_2)^{2k}t^{2k+1}} \int_0^{\frac{\beta}{c_1}} \bigl[c_2t+\beta-(c_1+c_2)t_1 \bigr]^{2k-1-j} \dif t_1+$$
$$ - \sum_{j=0}^{k-1} \binom{2k-1}{j}\frac{(2k+1)!}{(2k-1)!(2k-j)}\frac{(c_1t-\beta)^{2k-j}}{(c_1+c_2)^{2k}t^{2k+1}} \Bigl(\frac{c_2}{c_1}\Bigr)^{2k-1-2j}\int_0^{\frac{\beta}{c_1}} \bigl[c_2t+\beta-(c_1+c_2)t_1 \bigr]^{j} \dif t_1  =$$ 
$$ = \sum_{j=0}^{k-1} \binom{2k+1}{j+1}\frac{(c_1t-\beta)^{j+1}}{\bigl[(c_1+c_2)t\bigr]^{2k+1}}\Bigl[(c_2t +\beta)^{2k-j}- \Bigl(\frac{c_2}{c_1}\Bigr)^{2k-j}(c_1t-\beta)^{2k-j} \Bigr]+$$
$$ -\sum_{j=0}^{k-1} \binom{2k+1}{j+1}\frac{(c_1t-\beta)^{2k-j}}{\bigl[(c_1+c_2)t\bigr]^{2k+1}}\Bigl(\frac{c_2}{c_1}\Bigr)^{2k-1-2j}\Bigl[(c_2t +\beta)^{j+1}- \Bigl(\frac{c_2}{c_1}\Bigr)^{j+1}(c_1t-\beta)^{j+1} \Bigr]=$$
\begin{equation}\label{ultimoPassaggioDimMassimoRip+d}
= \sum_{j=1}^{k} \binom{2k+1}{j} \frac{ \Bigl[ (c_1t-\beta)^{j}(c_2t+\beta)^{2k+1-j}  - \Bigl(\frac{c_2}{c_1}\Bigr)^{2k+1-2j}(c_1t-\beta)^{2k+1-j}(c_2t+\beta)^{j} \Bigr]}{\bigl[(c_1+c_2)t \bigr]^{2k+1}}
\end{equation}
The term $j=0$ of (\ref{massimoRip+d}) is obtained by evaluating the second integral of (\ref{massimoRipIntegrale+d}) which yields
\begin{equation}
\int_0^{\frac{\beta}{c_1}}\dif t_1 \int_{\frac{c_1t-\beta}{c_1+c_2}+t_1}^t \frac{(2k+1)!}{(2k-1)!}\frac{(t-t_2)^{2k-1}}{t^{2k+1}} \dif t_2=\Bigl(\frac{c_2t+\beta}{(c_1+c_2)t}\Bigr)^{2k+1}-\Bigl(\frac{c_2}{c_1}\Bigr)^{2k+1}\Bigl(\frac{c_1t-\beta}{(c_1+c_2)t}\Bigr)^{2k+1}
\end{equation}
\hfill$\Box$

\begin{remark}
The distributions (\ref{massimoRip+p}) and (\ref{massimoRip+d}) can be unified in the following form, for $n\in \mathbb{N}, \beta \in [0, c_1t]$
\begin{equation}\label{massimoRip+n}
P\{\max_{0\le s\le t}\mathcal{T}(s)< \beta\ |\ N(t) = n,\ V(0) = c_1\} = 
\end{equation}
$$ = \frac{1}{\bigl[(c_1+c_2)t \bigr]^{n}} \sum_{j=0}^{\floor{\frac{n-1}{2}}} \binom{n}{j} \Bigl[ (c_1t-\beta)^{j}(c_2t+\beta)^{n-j}  - \Bigl(\frac{c_2}{c_1}\Bigr)^{n-2j}(c_1t-\beta)^{n-j}(c_2t+\beta)^{j} \Bigr]$$
Each term of (\ref{massimoRip+n}) is non-negative. Ideed, if we write the $j$-th term as
$$(c_1t-\beta)^j(c_2t+\beta)^j \Bigl[(c_2t+\beta)^{n-2j}-(c_2t -\frac{c_2}{c_1}\beta)^{n-2j} \Bigr]$$
and since $c_2t+\beta>c_2t -\frac{c_2}{c_1}\beta$ for $t, c_1,c_2>0$ and $\beta \in [0,c_1t]$ we conclude that (\ref{massimoRip+n}) is a sum of non-negative terms.
\hfill $\Diamond$
\end{remark}

\begin{remark}
For $c_1 = c_2 = c$ the distribution (\ref{massimoRip+n}) takes the form
\begin{equation}\label{massimoRip+nSimmetrico}
P\{\max_{0\le s\le t}\mathcal{T}(s)< \beta\ |\ N(t) = n,\ V(0) = c\} = \frac{1}{(2ct)^{n}} \sum_{j=0}^{\floor{\frac{n-1}{2}}} \binom{n}{j} (c^2t^2-\beta^2)^{j}\Bigl[(ct+\beta)^{n-2j}  - (ct-\beta)^{n-2j}\Bigr]
\end{equation}
The interested reader can prove that (\ref{massimoRip+nSimmetrico}) is the same for $n = 2k-1$ and $n = 2k$, $k \in \mathbb{N}$. This can be shown by considering the expression with $n = 2k-1$, by multiplying it for $\frac{(ct-\beta)+(ct+\beta)}{2ct}$ and by performing some calculations.
Formulas of Corollary $5.1$ of \cite{CO2020} confirm this cyclic behavior and furthermore they permit us to write down the following relationships, for integer $k\ge0,\ 0\le \beta \le ct$
\begin{equation}\label{massimoRip+pdSimmetricoRelazioni}
P\{\max_{0\le s\le t}\mathcal{T}(s)< \beta\ |\ V(0) = c,\ N(t) = 2k+1\} =  \int_0^\beta 2\frac{(2k+1)!}{k!^2} \frac{(c^2t^2-x^2)^{k}}{(2ct)^{2k+1}}\dif x =
\end{equation}
$$ = \frac{\beta}{ct}\sum_{j=0}^{k}\binom{2j}{j}\Bigl(\frac{c^2t^2-\beta^2}{c^2t^2} \Bigr)^{j} \frac{1}{2^{2j}} \ =\ \sum_{j=0}^{k} \binom{2k+1}{j} \frac{(c^2t^2-\beta^2)^{j}}{(2ct)^{2k+1}}\Bigl[(ct+\beta)^{2k+1-2j}  - (ct-\beta)^{2k+1-2j}\Bigr]  =$$
$$ = P\{\max_{0\le s\le t}\mathcal{T}(s)< \beta\ |\ N(t) = 2k+2,\ V(0) = c\} $$
In order to prove the third equality of (\ref{massimoRip+pdSimmetricoRelazioni}) we show the following unexpected relationship
\begin{equation}\label{massimoRip+pdSimmetricoRelazione2}
\int_0^\beta 2\frac{(2k+1)!}{k!^2} \frac{(c^2t^2-x^2)^{k}}{(2ct)^{2k+1}}\dif x= 
\end{equation}
$$ = \sum_{j=0}^{k} \binom{2k+1}{j} \frac{ \Bigl[ (ct-\beta)^{j}(ct+\beta)^{2k+1-j}  -(ct-\beta)^{2k+1-j}(ct+\beta)^{j} \Bigr]}{(2ct)^{2k+1}}$$
Equation (\ref{massimoRip+pdSimmetricoRelazione2}) can be proved by performing successive integration by parts. We obtain that
\begin{equation}\label{integraleNucleoDensita}
\int_0^\beta (ct-x)^{k}(ct+x)^{k}\dif x = \sum_{j=0}^k \frac{k!^2}{(k-j)!(k+j+1)!}\Bigl[(ct-\beta)^{k-j}(ct+\beta)^{k+j+1}-(ct)^{2k+1} \Bigr] =
\end{equation}
$$= -\sum_{j=0}^k \frac{k!^2}{(k-j)!(k+j+1)!}\Bigl[(ct-\beta)^{k+j+1}(ct+\beta)^{k-j}-(ct)^{2k+1} \Bigr]$$
where the second member is obtained by considering $(ct+x)$ as the integral term in the integration by parts and the third member is obtained by considering $(ct-x)$ as the integral term. Then
\begin{equation}\label{2integraleNucleoDensita}
2\int_0^\beta (ct-x)^{k}(ct+x)^{k}\dif x = 
\end{equation}
$$ = \sum_{j=0}^k \frac{k!^2}{(k-j)!(k+j+1)!}\Bigl[(ct-\beta)^{k-j}(ct+\beta)^{k+j+1}-(ct-\beta)^{k+j+1}(ct+\beta)^{k-j} \Bigr] =$$
\begin{equation}\label{2integraleNucleoDensita2}
= \sum_{j=0}^k \frac{k!^2}{j!(2k+1-j)!}\Bigl[(ct-\beta)^{j}(ct+\beta)^{2k+1-j}-(ct-\beta)^{2k+1-j}(ct+\beta)^{j} \Bigr]
\end{equation}
By multiplying (\ref{2integraleNucleoDensita}) and (\ref{2integraleNucleoDensita2}) by $ \frac{(2k+1)!}{k!^2}\frac{1}{(2ct)^{2k+1}}$ we obtain (\ref{massimoRip+pdSimmetricoRelazione2})
\hfill$\Diamond$
\end{remark}

\begin{remark}
In the frame of reference $(x',t')$ related to the original frame of reference $(x,t)$ by the relationships
\begin{equation}
\begin{cases}
x' = x+\frac{c_2-c_1}{2}t\\
t' = t
\end{cases}
\end{equation}
The distribution of the maximum (\ref{massimoRip+n}) (and thus the density) of the symmetric telegraph process coincides with (\ref{massimoRip+pdSimmetricoRelazione2}) with $c = \frac{c_1+c_2}{2},\ k=\floor{\frac{n-1}{2}}$. This can be proved by observing that
$$(c_1t-\beta) = \Bigl(\frac{c_1+c_2}{2}t'-\beta'\Bigr)\ ,\ \ \ \ \ (c_2t+\beta) = \Bigl(\frac{c_1+c_2}{2}t'+\beta'\Bigr)$$
Furthermore the velocities in the two frames of reference are related by 
$$\frac{\dif x'}{\dif t'}=\frac{\dif x}{\dif t}+\frac{c_2-c_1}{2}$$
Thus if in $(x,t),\ \frac{\dif x}{\dif t} = c_1$ (denoted below by $\frac{\dif x_+}{\dif t}$) in $(x',t')$ the observer sees the particle moving with velocity $\frac{\dif x_+'}{\dif t'}=\frac{c_1+c_2}{2}$. Clearly if $\frac{\dif x}{\dif t} = -c_2$ (denoted below by $\frac{\dif x_-}{\dif t}$) we have $\frac{\dif x_-'}{\dif t'}=-\frac{c_1+c_2}{2}$.
\\In (\ref{massimoRip+n}) we can also write
\begin{equation}
\Bigl( \frac{c_2}{c_1}\Bigr)^{n-2j} = \biggl(\frac{|\frac{\dif x_-}{\dif t}|}{\frac{\dif x_+}{\dif t} } \biggr)^{n-2j}  = \biggl(\frac{|\frac{\dif x_-'}{\dif t'}|}{\frac{\dif x_+'}{\dif t'} } \biggr)^{n-2j}  = 1
\end{equation}
The rightward point of the support in $(x',t')$ becomes $ \frac{c_1+c_2}{2}t'$.
\\In view of (\ref{massimoRip+pdSimmetricoRelazione2}), with $c =  \frac{c_1+c_2}{2}$ we obtain the distribution of the maximum in the symmetric case when $n = 2k+1$.
\hfill$\Diamond$
\end{remark}

\begin{remark}
We give the generating function of the cumulative distribution function (\ref{massimoRip+pdSimmetricoRelazioni}) as follows, for $|u|<1$
\begin{equation}
\sum_{k=0}^\infty u^kP\{\max_{0\le s\le t}\mathcal{T}(s)< \beta\ |\ V(0) = c,\ N(t) = 2k+1\} = \frac{\beta}{ct}\sum_{k=0}^\infty u^k\sum_{j=0}^{k}\binom{2j}{j}\Bigl(\frac{\sqrt{c^2t^2-\beta^2}}{2ct} \Bigr)^{2j}  = 
\end{equation}
$$ = \frac{\beta}{ct} \sum_{j=0}^\infty\binom{2j}{j}\Bigl(\frac{\sqrt{c^2t^2-\beta^2}}{2ct} \Bigr)^{2j}\sum_{k=j}^\infty u^k =\frac{\beta}{ct(1-u)} \sum_{j=0}^\infty\binom{2j}{j}\Bigl(\frac{\sqrt{u}\sqrt{c^2t^2-\beta^2}}{2ct} \Bigr)^{2j} = $$
\begin{equation}
=\frac{\beta}{(1-u)\sqrt{c^2t^2-u(c^2t^2-\beta^2)}}
\end{equation}
where in the last step we used the relationship
\begin{equation}
\frac{1}{\sqrt{1-x^2}} = \sum_{j=0}^\infty \frac{\Gamma\bigl(\frac{1}{2}\bigr)}{j!\Gamma\bigl(\frac{1}{2}-j\bigr)} \bigl(-x^2\bigr)^j = \sum_{j=0}^\infty \binom{2j}{j}\Bigl(\frac{x}{2}\Bigr)^{2j}
\end{equation}
for $|x|<1$.
\hfill$\Diamond$
\end{remark}

Interesting formulas can be given for the conditional density functions of the maximum.

\begin{corollary}\label{corollarioMassimo+p}
Let $\{\mathcal{T}(t)\}_{t\ge0}$ be the asymmetric telegraph process. For $k\in \mathbb{N},\ \beta\in (0, c_1t)$
\begin{equation}\label{massimo+p}
P\{\max_{0\le s\le t}\mathcal{T}(s)\in \dif \beta\ |\ V(0) = c_1,\ N(t) = 2k\}/ \dif \beta = 
\end{equation}
$$ = \frac{(2k)!}{k!(k-1)!}\frac{(c_1t-\beta)^{k-1}(c_2t+\beta)^{k} +\Bigl(\frac{c_2}{c_1}\Bigr)^{2} (c_1t-\beta)^{k}(c_2t+\beta)^{k-1} }{\bigl[(c_1+c_2)t \bigr]^{2k}} +$$
$$+\Bigl[ 1-\Bigl(\frac{c_1}{c_2}\Bigr)^2\ \Bigr]\sum_{j=0}^{k-2} \frac{(2k)!}{j!(2k-1-j)!}\Bigl(\frac{c_2}{c_1}\Bigr)^{2k-2j}\frac{(c_1t-\beta)^{2k-1-j}(c_2t+\beta)^{j}}{\bigl[(c_1+c_2)t \bigr]^{2k}}$$
\end{corollary}

\paragraph{Proof.}
The derivative of (\ref{massimoRip+p}) splits down into four terms as
\begin{equation}
\frac{\dif }{\dif \beta}P\{\max_{0\le s\le t}\mathcal{T}(s)<\beta\ |\ V(0) = c_1,\ N(t) = 2k\} = 
\end{equation}
$$= \frac{1}{\bigl[(c_1+c_2)t \bigr]^{2k}}  \sum_{j=0}^{k-1} \binom{2k}{j} \Bigl[ -j(c_1t-\beta)^{j-1}(c_2t+\beta)^{2k-j}  +(2k-j)(c_1t-\beta)^{j}(c_2t+\beta)^{2k-j-1} +$$
$$+(2k-j)\Bigl(\frac{c_2}{c_1}\Bigr)^{2k-2j}(c_1t-\beta)^{2k-j-1}(c_2t+\beta)^{j} -j\Bigl(\frac{c_2}{c_1}\Bigr)^{2k-2j}(c_1t-\beta)^{2k-j}(c_2t+\beta)^{j-1} \Bigr]=$$
$$ = \frac{1}{\bigl[(c_1+c_2)t \bigr]^{2k}} \Biggl[-\sum_{j=0}^{k-2}\frac{(2k)!}{j!(2k-1-j)!}(c_1t-\beta)^{j}(c_2t+\beta)^{2k-j-1}  +  $$
$$+\sum_{j=0}^{k-1}\frac{(2k)!}{j!(2k-1-j)!}(c_1t-\beta)^{j}(c_2t+\beta)^{2k-j-1} +\sum_{j=0}^{k-1}\frac{(2k)!}{j!(2k-1-j)!}\Bigl(\frac{c_2}{c_1}\Bigr)^{2k-2j}(c_1t-\beta)^{2k-j-1}(c_2t+\beta)^{j}+$$
$$-\sum_{j=0}^{k-2} \frac{(2k)!}{j!(2k-1-j)!}\Bigl(\frac{c_2}{c_1}\Bigr)^{2k-2j-2}(c_1t-\beta)^{2k-j-1}(c_2t+\beta)^{j} \Biggr] =  $$
$$ =\frac{1}{\bigl[(c_1+c_2)t \bigr]^{2k}} \Biggl[ \frac{(2k)!}{k!(k-1)!}(c_1t-\beta)^{k-1}(c_2t+\beta)^{k} +\frac{(2k)!}{k!(k-1)!}\Bigl(\frac{c_2}{c_1}\Bigr)^{2} (c_1t-\beta)^{k}(c_2t+\beta)^{k-1} +$$
$$+\Bigl[ 1-\Bigl(\frac{c_1}{c_2}\Bigr)^2\ \Bigr]\sum_{j=0}^{k-2} \frac{(2k)!}{j!(2k-1-j)!}\Bigl(\frac{c_2}{c_1}\Bigr)^{2k-2j}(c_1t-\beta)^{2k-1-j}(c_2t+\beta)^{j}\Biggr] $$
\hfill$\Box$
\\

For $c_1 = c_2 = c$ we have
\begin{equation}\label{massimo+pSimmetrico}
P\{\max_{0\le s\le t}\mathcal{T}(s)\in \dif \beta\ |\ V(0) = c,\ N(t) = 2k\} = 2\frac{(2k-1)!}{(k-1)!^2} \frac{(c^2t^2-\beta^2)^{k-1}}{(2ct)^{2k-1}}\dif \beta =
\end{equation}
$$ = 2P\{\mathcal{T}(t)\in \dif \beta\ |\ N(t) = 2k\} $$
which coincides with formula ($3.10$) of \cite{CO2020} with $k+1$ replaced by $k$.
\\

For the case where $V(0) = c_1,\ N(t) = 2k+1$, the density of the maximum of the asymmetric telegraph process is presented in the next corollary.

\begin{corollary}\label{corollarioMassimo+d}
Let $\{\mathcal{T}(t)\}_{t\ge0}$ be the asymmetric telegraph process. For $k\in \mathbb{N}_0,\ \beta\in (0, c_1t)$
\begin{equation}\label{massimo+d}
P\{\max_{0\le s\le t}\mathcal{T}(s)\in \dif \beta\ |\ V(0) = c_1,\ N(t) = 2k+1\}/ \dif \beta = 
\end{equation}
$$ =\ \Bigl(1+\frac{c_2}{c_1} \Bigr)\frac{(2k+1)!}{k!^2}\frac{(c_1t-\beta)^{k}(c_2t+\beta)^{k} }{\bigl[(c_1+c_2)t \bigr]^{2k+1}} \ +$$
$$+\ \Bigl[ 1-\Bigl(\frac{c_1}{c_2}\Bigr)^2\ \Bigr]\sum_{j=0}^{k-1} \frac{(2k+1)!}{j!(2k-j)!}\Bigl(\frac{c_2}{c_1}\Bigr)^{2k+1-2j}\frac{(c_1t-\beta)^{2k-j}(c_2t+\beta)^{j}}{\bigl[(c_1+c_2)t \bigr]^{2k+1}}$$
\end{corollary}

\paragraph{Proof.}
The proof is similar to the proof of Corollary \ref{corollarioMassimo+p} and therefore it is omitted.
\hfill$\Box$
\\

In the symmetric case, $c_1 = c_2 = c$,
\begin{equation}\label{massimo+dSimmetrico}
P\{\max_{0\le s\le t}\mathcal{T}(s)\in \dif \beta\ |\ V(0) = c,\ N(t) = 2k+1\} = 2\frac{(2k+1)!}{k!^2} \frac{(c^2t^2-\beta^2)^{k}}{(2ct)^{2k+1}}\dif \beta =
\end{equation}
$$ = 2P\{\mathcal{T}(t)\in \dif \beta\ |\ N(t) = 2k+1\} $$
By comparing (\ref{massimo+pSimmetrico}) and (\ref{massimo+dSimmetrico}) we have that, for $k\in \mathbb{N},\ 0\le \beta\le ct$
\begin{equation}\label{catenaEquazioniMassimo+}
P\{\max_{0\le s\le t}\mathcal{T}(s)\in \dif \beta | V(0) = c, N(t) = 2k-1\} =P\{\max_{0\le s\le t}\mathcal{T}(s)\in \dif \beta | V(0) = c, N(t) = 2k\} =
\end{equation}
$$ = 2P\{\mathcal{T}(t)\in \dif \beta\ |\ N(t) = 2k\}  =2P\{\mathcal{T}(t)\in \dif \beta\ |\ N(t) = 2k-1\} $$
The equations in (\ref{catenaEquazioniMassimo+}) describe two incredible facts about the symmetric telegraph process: when the motion starts with positive initial velocity the conditional distributions of the maximum have a cyclic behavior and some kind of reflection principle holds.
\\Analogous regularities do not hold in the asymmetric case.

\begin{remark}
The distribution function of the maximum when $N(t)=n$ is a polynomial function of order $n$ while the associated density function is of order $n-1$. The coefficients depend on the possible velocities $c_1$ and $c_2$ in a rather complicated way.
For $N(t)=3,\ V(0)=c_1$, for example we have that
\begin{equation}
P\{\max_{0\le s \le t}\mathcal{T}(t)<\beta\ |\ V(0) = c_1,\ N(t) = 3\}= \frac{\beta^3}{(c_1t)^3}+\frac{3\beta(c_1t-\beta)(c_2t+\beta)}{(c_1+c_2)\,c_1^2\,t^3}
\end{equation}
for $0\le \beta \le c_1t$, and
\begin{equation}\label{densita3+}
P\{\max_{0\le s \le t}\mathcal{T}(t)\in \dif\beta\ |\ V(0) = c_1,\ N(t) = 3\}/ \dif \beta\ =\ \frac{3(c_1t-\beta)(c_2t+\beta)}{(c_1+c_2)\,c_1^2\,t^3}+\frac{3\beta (c_1-c_2)(c_1t-\beta)}{(c_1+c_2)(c_1t)^3}
\end{equation}
for $0<\beta<c_1t$. The density (\ref{densita3+}) displays a maximal point at 
\begin{equation}\label{puntoDiMassimo+3}
\beta_{max} = \frac{c_1(c_1-c_2)t}{2c_1-c_2} \ \ \ \ \ \ \ if\ \ c_1>c_2 \ \ or \ \ c_1<\frac{c_2}{2}
\end{equation}
otherwise it is monotonically decreasing from $\beta=0$, with value $\frac{3c_2}{(c_1+c_2)c_1t}$, to $\beta = c_1t$ where it vanishes. In the maximal point $\beta_{max}$ in (\ref{puntoDiMassimo+3}) the density function takes the value
$$\frac{3c_1}{(c_1+c_2)(2c_1-c_2)t} = \frac{3}{(c_1^2-c_2^2)\,t^2}\beta_{max}$$
We recall that the study of the distributions for $k$ small, both in the odd and the even case, has been crucial to identify the general form of the cumulative distribution we have presented in Theorem \ref{teoremaMassimoRip+p} and Theorem \ref{teoremaMassimoRip+d}.
\hfill $\Diamond$
\end{remark}

We conclude this section with the unconditional distribution of the maximum of the telegraph process starting with velocity $c_1>0$. In this case there is a positive probability mass in the point $\beta = c_1t$ where the moving particle arrives if no Poisson event occurs in the time interval $[0,t]$.
\begin{theorem}\label{teoremaMassimoRip+}
Let $\{\mathcal{T}(t)\}_{t\ge0}$ be the asymmetric telegraph process. For $\beta\in [0, c_1t]$ we have that
\begin{equation}\label{massimoRip+}
P\{\max_{0\le s\le t}\mathcal{T}(s)< \beta\ |\ V(0) = c_1\} = 
\end{equation}
$$ =e^{-\lambda t}\sum_{r=1}^\infty I_r\Bigl(\frac{2\lambda}{c_1+c_2} \sqrt{(c_1t-\beta)(c_2t+\beta)}\Bigr)\Biggl[ \Biggl(\sqrt{\frac{c_2t+\beta}{c_1t-\beta}}\Biggr)^r-\Biggl(\frac{c_2}{c_1}\sqrt{\frac{c_1t-\beta}{c_2t+\beta}}\Biggr)^r\ \Biggr]$$
and
\begin{equation}\label{massimoRip+Singolarita}
P\{\max_{0\le s\le t}\mathcal{T}(s)=c_1t\ |\ V(0) = c_1\} =e^{-\lambda t}
\end{equation}
\end{theorem}

\paragraph{Proof.}
In view of (\ref{massimoRip+p}) we evaluate the following joint probability
\begin{equation}\label{massimoRip+UnitoPari}
P\{\max_{0\le s\le t}\mathcal{T}(s)< \beta,\ \bigcup_{k=1}^\infty \{N(t) = 2k\}\ |\ V(0) = c_1\} = 
\end{equation}
$$=\sum_{k = 1}^{\infty}e^{-\lambda t}\frac{(\lambda t)^{2k}}{(2k)!}  \sum_{j=0}^{k-1} \binom{2k}{j} \, \frac{(c_1t-\beta)^{j}(c_2t+\beta)^{2k-j}  - \Bigl(\frac{c_2}{c_1}\Bigr)^{2k-2j}(c_1t-\beta)^{2k-j}(c_2t+\beta)^{j}}{\bigl[(c_1+c_2)t \bigr]^{2k}}= $$
$$ = e^{-\lambda t}  \sum_{j=0}^\infty \frac{1}{j!} \sum_{k = j+1}^{\infty} \Bigl( \frac{\lambda}{c_1+c_2}\Bigr)^{2k}\, \frac{(c_1t-\beta)^{j}(c_2t+\beta)^{2k-j}  - \Bigl(\frac{c_2}{c_1}\Bigr)^{2k-2j}(c_1t-\beta)^{2k-j}(c_2t+\beta)^{j}}{(2k-j)!} =$$
$$ = e^{-\lambda t}  \sum_{j=0}^\infty \frac{1}{j!} \sum_{k = 0}^{\infty} \Bigl( \frac{\lambda}{c_1+c_2}\Bigr)^{2k+2j+2}\, \frac{1}{(2k+2+j)!} \Bigl[(c_1t-\beta)^{j}(c_2t+\beta)^{2k+2+j} + $$
$$-\ \Bigl(\frac{c_2}{c_1}\Bigr)^{2k+2}(c_1t-\beta)^{2k+2+j}(c_2t+\beta)^{j}\Bigr]\ =$$
\begin{equation} =e^{-\lambda t}\sum_{k=0}^\infty I_{2k+2}\Bigl(\frac{2\lambda}{c_1+c_2} \sqrt{(c_1t-\beta)(c_2t+\beta)}\Bigr)\Biggl[ \Biggl(\sqrt{\frac{c_2t+\beta}{c_1t-\beta}}\Biggr)^{2k+2}-\Biggl(\frac{c_2}{c_1}\sqrt{\frac{c_1t-\beta}{c_2t+\beta}}\Biggr)^{2k+2}\: \Biggr]
\end{equation}
In light of (\ref{massimoRip+d}), with similar steps we have that
\begin{equation}\label{massimoRip+UnitoDispari}
P\{\max_{0\le s\le t}\mathcal{T}(s)< \beta,\ \bigcup_{k=0}^\infty\{ N(t) = 2k+1\}\ |\ V(0) = c_1\} = 
\end{equation}
$$ =e^{-\lambda t}\sum_{k=0}^\infty I_{2k+1}\Bigl(\frac{2\lambda}{c_1+c_2} \sqrt{(c_1t-\beta)(c_2t+\beta)}\Bigr)\Biggl[ \Biggl(\sqrt{\frac{c_2t+\beta}{c_1t-\beta}}\Biggr)^{2k+1}-\Biggl(\frac{c_2}{c_1}\sqrt{\frac{c_1t-\beta}{c_2t+\beta}}\Biggr)^{2k+1}\ \Biggr]$$
The sum of (\ref{massimoRip+UnitoPari}) and (\ref{massimoRip+UnitoDispari}) yields the claimed result (\ref{massimoRip+}).
\\It can be checked that
\begin{equation}
P\{\max_{0\le s\le t}\mathcal{T}(s)< c_1t\ |\ V(0) = c_1\} = e^{-\lambda t} (e^{\lambda t}-1) = 1- e^{-\lambda t} = 
\end{equation}
$$ = 1- P\{\max_{0\le s\le t}\mathcal{T}(s)=c_1t\ |\ V(0) = c_1\}$$
\hfill$\Box$

\section{Asymmetric telegraph process with negative initial velocity}
The results of the previous section permit us to derive the explicit distribution of the maximum of the asymmetric telegraph process when the initial velocity is negative.
\\In this case the probability mass at $\beta=0$, i.e.
\begin{equation}\label{InizioCap3massimo-pSingolarita}
P\{\max_{0\le s\le t}\mathcal{T}(s)=0\ |\ V(0) = -c_2,\ N(t) = n\} 
\end{equation}
with $n\in \mathbb{N}$, is strictly positive and independent of $t$. The probabilites (\ref{InizioCap3massimo-pSingolarita}) display a ciclicity which is not valid for the corresponding cumulative distribution functions. This means that, for all natural numbers $k \ge 1$ we have that
$$P\{\max_{0\le s\le t}\mathcal{T}(s)=0 | V(0) = -c_2, N(t) = 2k-1\} = P\{\max_{0\le s\le t}\mathcal{T}(s)=0 | V(0) = -c_2, N(t) = 2k\}  $$
whereas
$$P\{\max_{0\le s\le t}\mathcal{T}(s)\le \beta | V(0) = -c_2, N(t) = 2k-1\}\not = P\{\max_{0\le s\le t}\mathcal{T}(s)\le \beta | V(0) = -c_2, N(t) = 2k\}  $$

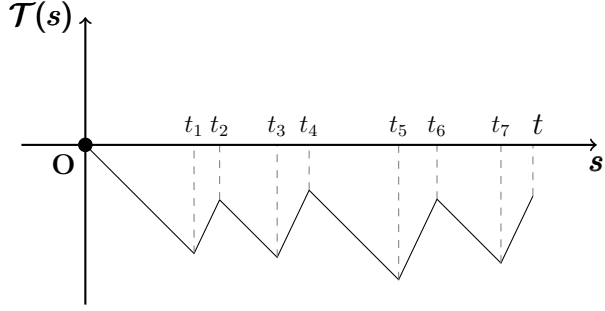
\begin{figure}
\begin{minipage}{1\textwidth}
\centering
\begin{tikzpicture}[scale = 0.85]
\draw(0,0) --(1.7, -1.7)--(2.1, -0.86)--(3, -1.76)--(3.5, -0.71)--(4.9, -2.11)--(5.5, -0.85)--(6.5, -1.85)--(7, -0.8);;
\draw[dashed, gray](1.7, -1.7)-- (1.7,0) node[above, black, scale = 0.8]{$t_1$};
\draw[dashed, gray] (2.1, -0.86)-- (2.1,0) node[above, black, scale = 0.8]{$t_2$};
\draw[dashed, gray] (3, -1.76)-- (3,0) node[above, black, scale = 0.8]{$t_3$};
\draw[dashed, gray] (3.5, -0.71)-- (3.5,0) node[above, black, scale = 0.8]{$t_4$};
\draw[dashed, gray] (4.9, -2.11)-- (4.9,0) node[above, black, scale = 0.8]{$t_5$};
\draw[dashed, gray] (5.5, -0.85)-- (5.5,0) node[above, black, scale = 0.8]{$t_6$};
\draw[dashed, gray] (6.5, -1.85)-- (6.5,0) node[above, black, scale = 0.8]{$t_7$};
\draw[dashed, gray] (7,-0.8)-- (7,0) node[above, black, scale = 1]{$\ t$};
\draw[->, thick] (-1,0) -- (8,0) node[below, scale = 1]{$\pmb{s}$};
\draw[->, thick] (0,-2.5) -- (0,2) node[left, scale = 1]{ $\pmb{\mathcal{T}(s)}$};
\filldraw[black] (0,0) circle (3pt) node[below left, black, scale =0.9]{\textbf{O}};
\end{tikzpicture}
\caption{motion with negative initial velocity, $c_1>c_2$}\label{traiettoria0}
\end{minipage}
\end{figure}

We start our analysis with the following result.

\begin{theorem}\label{teoremaMassimoRip-p}
Let $\{\mathcal{T}(t)\}_{t\ge0}$ be the asymmetric telegraph process. For $k\in \mathbb{N}_0,\ \beta\in [0, c_1t]$
\begin{equation}\label{massimoRip-p}
P\{\max_{0\le s\le t}\mathcal{T}(s)\le \beta\ |\ V(0) = -c_2,\ N(t) = 2k\} = 
\end{equation}
$$ =\ \sum_{j=0}^{k} \binom{2k}{j} \frac{ (c_1t-\beta)^{j}(c_2t+\beta)^{2k-j}}{\bigl[(c_1+c_2)t \bigr]^{2k}} \  -\ \sum_{j=0}^{k-1} \binom{2k}{j}\Bigl(\frac{c_2}{c_1}\Bigr)^{2k-1-2j}\frac{ (c_1t-\beta)^{2k-j}(c_2t+\beta)^{j} }{\bigl[(c_1+c_2)t \bigr]^{2k}}$$
\end{theorem}

\paragraph{Proof.}
At time $T_1 = t_1$ we have two possible situations:
\begin{enumerate}
\item $\beta-(-c_2t_1)\ge c_1(t-t_1)$ ;
\item $\beta-(-c_2t_1)< c_1(t-t_1)$ .
\end{enumerate}
In the case 1 the particle went so far to the left that it will not be able to cross the level $\beta$ in the remaining time interval $(t_1,t)$.
\\In the case 2 the particle can, with positive probability, overcome the level $\beta$ during the time interval $(t_1,t)$.
\\Thus, by considering also the distribution (\ref{massimoRip+d}) and arguing by induction as above, we can write 
\begin{equation}
P\{\max_{0\le s\le t}\mathcal{T}(s)\le \beta\ |\ V(0) = -c_2,\ N(t) = 2k\} = 
\end{equation}
$$ = \int_0^{\frac{c_1t-\beta}{c_1+c_2}} P\{\max_{0\le s\le t-t_1}\mathcal{T}(s)\le \beta+c_2t_1\ |\ V(0) = c_1,\ N(t-t_1) = 2k-1\}\ \cdot$$
$$\cdot\ P\{T_1\in \dif t_1\ |\ N(t) = 2k\}\ +\ \int_{\frac{c_1t-\beta}{c_1+c_2}}^t P\{T_1\in \dif t_1\ |\ N(t) = 2k\}=$$
$$ =\  \int_0^{\frac{c_1t-\beta}{c_1+c_2}} \frac{2k(t-t_1)^{2k-1}}{t^{2k}}\frac{\dif t_1}{\bigl[(c_1+c_2)(t-t_1) \bigr]^{2k-1}}\cdot$$
$$ \cdot\sum_{j=0}^{k-1} \binom{2k-1}{j} \Biggl[ \bigl[c_1(t-t_1)-(\beta+c_2t_1)\bigr]^{j}\bigl[c_2(t-t_1)+\beta+c_2t_1\bigr]^{2k-1-j}\  +$$
$$-\ \Bigl(\frac{c_2}{c_1}\Bigr)^{2k-1-2j}\bigl[c_1(t-t_1)-(\beta+c_2t_1)\bigr]^{2k-1-j}\bigl[c_2(t-t_1)+\beta+c_2t_1\bigr]^{j} \Biggr]\ +\ \Bigl(\frac{c_2t+\beta}{(c_1+c_2)t}\Bigr)^{2k}=$$
$$ =\ \sum_{j=0}^{k-1} 2k\binom{2k-1}{j} \frac{1}{(c_1+c_2)^{2k-1}t^{2k}} \Biggl[ (c_2t+\beta)^{2k-1-j}  \int_0^{\frac{c_1t-\beta}{c_1+c_2}} \bigl[c_1t-\beta-(c_1+c_2)t_1\bigr]^{j}\dif t_1\ +$$
$$-\ \Bigl(\frac{c_2}{c_1}\Bigr)^{2k-1-2j}(c_2t+\beta)^{j}\int_0^{\frac{c_1t-\beta}{c_1+c_2}} \bigl[c_1t-\beta-(c_1+c_2)t_1\bigr]^{2k-1-j}\dif t_1\Biggr]+ \Bigl(\frac{c_2t+\beta}{(c_1+c_2)t}\Bigr)^{2k}$$
Some further steps yield (\ref{massimoRip-p}).
\hfill$\Box$

\begin{remark}\label{remarkSingolarita-p}
From (\ref{massimoRip-p}) we obtain the probability mass in $\beta=0$. For $k\in\mathbb{N}_0$
\begin{equation}\label{massimo-pSingolarita1}
P\{\max_{0\le s\le t}\mathcal{T}(s)=0\ |\ V(0) = -c_2,\ N(t) = 2k\} \ =\ \binom{2k}{k}\frac{(c_1c_2)^k}{(c_1+c_2)^{2k}}+ \Bigl(1-\frac{c_1}{c_2}\Bigr)\sum_{j=0}^{k-1} \binom{2k}{j} \frac{c_1^{j}\,c_2^{2k-j}}{(c_1+c_2)^{2k}}
\end{equation}
Alternatively, we can obtain
$$ P\{\max_{0\le s\le t}\mathcal{T}(s)=0\ |\ V(0) = -c_2,\ N(t) = 2k\}\ =\ \sum_{j=0}^{k} \binom{2k}{j} \frac{c_1^{j}\,c_2^{2k-j}}{(c_1+c_2)^{2k}}-\sum_{j=1}^{k} \binom{2k}{j-1} \frac{c_1^{j}\,c_2^{2k-j}}{(c_1+c_2)^{2k}}=$$
\begin{equation}\label{massimo-pSingolarita2}
 = \sum_{j=0}^{k}\Biggl[ \binom{2k}{j}-\binom{2k}{j-1}  \Biggr]\frac{c_1^{j}\,c_2^{2k-j}}{(c_1+c_2)^{2k}}
\end{equation}
We observe that the probability of the singular point depends only on the ratio of the two possible velocities. Let $c_1 = \alpha c_2,\ \alpha>0$, then we obtain
\begin{equation}\label{massimo-pSingolaritaAlpha}
P\{\max_{0\le s\le t}\mathcal{T}(s)=0\ |\ V(0) = -c_2,\ N(t) = 2k\} \ =\ \binom{2k}{k}\frac{\alpha^k}{(1+\alpha)^{2k}}+ (1-\alpha)\sum_{j=0}^{k-1} \binom{2k}{j} \frac{\alpha^{j}}{(1+\alpha)^{2k}}
\end{equation}
which is positive for all $\alpha>0$.
\\Finally, for $c_1 = c_2 = c$ (i.e. $\alpha = 1$) we have that
\begin{equation}
P\{\max_{0\le s\le t}\mathcal{T}(s)=0\ |\ N(t) = 2k,\ V(0) = -c\} = \binom{2k}{k}\frac{1}{2^{2k}}
\end{equation}
\hfill$\Diamond$
\end{remark}

\begin{theorem}\label{teoremaMassimoRip-d}
Let $\{\mathcal{T}(t)\}_{t\ge0}$ be the asymmetric telegraph process. For $k\in \mathbb{N}_0,\ \beta\in [0, c_1t]$
\begin{equation}\label{massimoRip-d}
P\{\max_{0\le s\le t}\mathcal{T}(s)\le \beta\ |\ V(0) = -c_2,\ N(t) = 2k+1\} = 
\end{equation}
$$ = \sum_{j=0}^{k} \binom{2k+1}{j} \frac{ (c_1t-\beta)^{j}(c_2t+\beta)^{2k+1-j}}{\bigl[(c_1+c_2)t \bigr]^{2k+1}} \  -\ \sum_{j=0}^{k-1} \binom{2k+1}{j}\Bigl(\frac{c_2}{c_1}\Bigr)^{2k-2j}\frac{ (c_1t-\beta)^{2k+1-j}(c_2t+\beta)^{j} }{\bigl[(c_1+c_2)t \bigr]^{2k+1}}$$
\end{theorem}

\paragraph{Proof.}
By writing
\begin{equation}
P\{\max_{0\le s\le t}\mathcal{T}(s)\le \beta\ |\ V(0) = -c_2,\ N(t) = 2k+1\} = 
\end{equation}
$$ = \int_0^{\frac{c_1t-\beta}{c_1+c_2}} P\{\max_{0\le s\le t-t_1}\mathcal{T}(s)\le \beta+c_2t_1\ |\ V(0) = c_1,\ N(t-t_1) = 2k\}\ \cdot $$
$$ \cdot \ P\{T_1\in \dif t_1\ |\ N(t) = 2k+1\}\ +\ \int_{\frac{c_1t-\beta}{c_1+c_2}}^t P\{T_1\in \dif t_1\ |\ N(t) = 2k+1\}$$
and by applying formula (\ref{massimoRip+p}), we obtain result (\ref{massimoRip-d}) by means of the same steps as in Theorem \ref{teoremaMassimoRip-p}.
\hfill$\Box$

\begin{remark}
By means of considerations similar to those of Remark \ref{remarkSingolarita-p}, for $\beta = 0, \ k\in \mathbb{N}_0$, we have that
$$ P\{\max_{0\le s\le t}\mathcal{T}(s)=0\ |\ V(0) = -c_2,\ N(t) = 2k+1\} =$$
\begin{equation}\label{massimo-dSingolarita1}
= \binom{2k+1}{k}\frac{c_1^k\,c_2^{k+1}}{(c_1+c_2)^{2k+1}}+ \Bigl(1-\frac{c_1}{c_2}\Bigr)\sum_{j=0}^{k-1} \binom{2k+1}{j} \frac{c_1^{j}\,c_2^{2k+1-j}}{(c_1+c_2)^{2k+1}}=
\end{equation}
\begin{equation}\label{massimo-dSingolarita2}
 = \sum_{j=0}^{k}\Biggl[ \binom{2k+1}{j}-\binom{2k+1}{j-1}  \Biggr]\frac{c_1^{j}\,c_2^{2k+1-j}}{(c_1+c_2)^{2k+1}}
\end{equation}
\hfill $\Diamond$
\end{remark}

We now show that the probability mass at $\beta = 0$ of the asymmetric telegraph process has a cyclic behavior.

\begin{proposition}
Let $\{\mathcal{T}(t)\}_{t\ge0}$ be the asymmetric telegraph process. For $k\in \mathbb{N}$
\begin{equation}\label{relazioneSingolarita}
P\{\max_{0\le s\le t}\mathcal{T}(s)=0\ |\ V(0) = -c_2,\ N(t) = 2k-1\} =
\end{equation}
$$ = P\{\max_{0\le s\le t}\mathcal{T}(s)=0\ |\ V(0) = -c_2,\ N(t) = 2k\}$$
\end{proposition}

\paragraph{Proof.}
In view of (\ref{massimo-dSingolarita2}) and (\ref{massimo-pSingolarita2}) the statement (\ref{relazioneSingolarita}) is equivalent to
\begin{equation}\label{uguaglianzaMassimoSingolarita}
\sum_{j=0}^{k-1}\Biggl[ \binom{2k-1}{j}-\binom{2k-1}{j-1}  \Biggr]\frac{c_1^{j}\,c_2^{2k-1-j}}{(c_1+c_2)^{2k-1}}=\sum_{j=0}^{k}\Biggl[ \binom{2k}{j}-\binom{2k}{j-1}  \Biggr]\frac{c_1^{j}\,c_2^{2k-j}}{(c_1+c_2)^{2k}}
\end{equation}
The first term of (\ref{uguaglianzaMassimoSingolarita}) can be written as
$$\sum_{j=0}^{k-1}\Biggl[ \binom{2k-1}{j}-\binom{2k-1}{j-1}  \Biggr]\frac{c_1^{j}\,c_2^{2k-1-j}}{(c_1+c_2)^{2k}}(c_2+c_1) =$$
$$= \sum_{j=0}^{k-1}\Biggl[ \binom{2k-1}{j}-\binom{2k-1}{j-1}  \Biggr]\frac{c_1^{j}\,c_2^{2k-j}}{(c_1+c_2)^{2k}} \:+\: \sum_{j=1}^{k}\Biggl[ \binom{2k-1}{j-1}-\binom{2k-1}{j-2}  \Biggr]\frac{c_1^{j}\,c_2^{2k-j}}{(c_1+c_2)^{2k}} =$$
$$ = \sum_{j=1}^{k-1}\Biggl[ \binom{2k-1}{j}-\binom{2k-1}{j-2}  \Biggr]\frac{c_1^{j}\,c_2^{2k-j}}{(c_1+c_2)^{2k}} \:+\: \Bigl(\frac{c_2}{c_1+c_2}\Bigr)^{2k} \: +\:\Biggl[ \binom{2k-1}{k-1}-\binom{2k-1}{k-2}  \Biggr]\frac{c_1^k\,c_2^k}{(c_1+c_2)^{2k}}  = $$
$$=\ \sum_{j=0}^{k}\Biggl[ \binom{2k}{j}-\binom{2k}{j-1}  \Biggr]\frac{c_1^{j}\,c_2^{2k-j}}{(c_1+c_2)^{2k}}$$
where in the last step we applied the following relationships
$$\binom{2k-1}{j}-\binom{2k-1}{j-2}  = \binom{2k}{j}-\binom{2k}{j-1} $$
and
$$ \binom{2k-1}{k-1}-\binom{2k-1}{k-2} =\binom{2k}{k}-\binom{2k}{k-1}$$
for positive integers $1\le j\le k-1$.
\hfill$\Box$
\\

We note that (\ref{relazioneSingolarita}) is the unique case where for the asymmetric telegraph process the cyclic behavior of the distribution of the maximum holds.

\begin{remark}
The probability (\ref{massimo-dSingolarita2}) can also be written as
\begin{equation}\label{massimo-SingolaritaA}
P\{\max_{0\le s\le t}\mathcal{T}(s)=0\ |\ V(0) = -c_2,\ N(t) = 2k+1\} = \Bigl(\frac{c_2}{c_1+c_2} \Bigr)^{k+1}\sum_{j=0}^k A_j^{(k)} \Bigl(\frac{c_1}{c_1+c_2} \Bigr)^{j}
\end{equation}
where the numbers $A_j^{(k)}$ form the following triangular matrix.
\begin{center}
\begin{tabular}{ c|c c c c c c c}
$\ \ k\ \backslash \ j$  & 0  & 1 & 2 & 3 & 4 & 5 & - \\
\hline
0 & 1 \\
1 & 1 & 1 \\
2 & 1 & 2 & 2 \\
3 & 1 & 3 & 5 & 5 \\
4 & 1 & 4 & 9 & 14 & 14 \\
5 & 1 & 5 & 14 & 28 & 42 & 42 \\
- & - & - & - & - & - & - & -  \\
\end{tabular}
\\
\centering
\small{Table 1. Triangular matrix of the coefficients $A_j^{(k)}, \ k\ge j\in \mathbb{N}_0$}
\end{center}
The numbers $A_j^{(k)}$ are related among themselves by the following recurrence relationships
\begin{equation}
A_0^{(0)}=1,\ \ \ A_k^{(k)} = A_{k-1}^{(k)},\ \ \ A_j^{(k)} = \sum_{i=0}^j A_i^{(k-1)},\ \ \ \ k>j\ge0
\end{equation}
\hfill $\Diamond$
\end{remark}

\begin{corollary}\label{corollarioMassimo-p}
Let $\{\mathcal{T}(t)\}_{t\ge0}$ be the asymmetric telegraph process. For $k\in \mathbb{N},\ \beta\in (0, c_1t)$
\begin{equation}\label{massimo-p}
P\{\max_{0\le s\le t}\mathcal{T}(s)\in \dif \beta\ |\ V(0) = -c_2,\ N(t) = 2k\}/ \dif \beta = 
\end{equation}
$$ =\ \Bigl(1+\frac{c_2}{c_1} \Bigr)\frac{(2k)!}{k!(k-1)!}\frac{(c_1t-\beta)^{k}(c_2t+\beta)^{k-1} }{\bigl[(c_1+c_2)t \bigr]^{2k}}\ +$$
$$+\ \Bigl[ 1-\Bigl(\frac{c_1}{c_2}\Bigr)^2\ \Bigr]\sum_{j=0}^{k-2} \frac{(2k)!}{j!(2k-1-j)!}\Bigl(\frac{c_2}{c_1}\Bigr)^{2k-1-2j}\frac{(c_1t-\beta)^{2k-1-j}(c_2t+\beta)^{j}}{\bigl[(c_1+c_2)t \bigr]^{2k}}$$
\end{corollary}

\paragraph{Proof.}
By working on formula (\ref{massimoRip-p}) and performing steps similar to those of Corollary \ref{corollarioMassimo+p} we obtain result (\ref{massimo-p}).
\hfill$\Box$
\\

We observe that for $c_1 = c_2 = c$ the sum appearing in (\ref{massimo-p}) is canceled out and we obtain that
\begin{equation}\label{massimoSimmetrico-p}
P\{\max_{0\le s\le t}\mathcal{T}(s)\in \dif \beta\ |\ N(t) = 2k,\ V(0) = -c\} = 2\frac{(2k)!}{k!(k-1)!}\frac{(c^2t^2-\beta^2)^{k-1}}{(2ct)^{2k}}(ct-\beta) \dif \beta = 
\end{equation}
$$ =2P\{\mathcal{T}(t)\in \dif \beta\ |\ N(t) = 2k,\ V(0) = -c\} $$
Formula (\ref{massimoSimmetrico-p}) shows that in the symmetric motion, if an even number of reversals occur, then a reflection principle holds for the absolutely continuous component in (0,ct). In the asymmetric case, the different absolute values of the velocities of motion, $c_1\not =c_2$, do not allow a similar property.

\begin{corollary}\label{corollarioMassimo-d}
Let $\{\mathcal{T}(t)\}_{t\ge0}$ be the asymmetric telegraph process. For $k\in \mathbb{N}_0,\ \beta\in (0, c_1t)$
\begin{equation}\label{massimo-d}
P\{\max_{0\le s\le t}\mathcal{T}(s)\in \dif \beta\ |\ V(0) = -c_2,\ N(t) = 2k+1\}/ \dif \beta = 
\end{equation}
$$ =\ \frac{(2k+1)!}{\bigl[(c_1+c_2)t \bigr]^{2k+1}}\Biggl[\frac{(c_1t-\beta)^{k}(c_2t+\beta)^{k} }{k!^2} +\Bigl(\frac{c_2}{c_1}\Bigr)^{2} \frac{(c_1t-\beta)^{k+1}(c_2t+\beta)^{k-1} }{(k+1)!(k-1)!} \ +$$
$$+\ \Bigl[ 1-\Bigl(\frac{c_1}{c_2}\Bigr)^2\Bigr]\sum_{j=0}^{k-2} \Bigl(\frac{c_2}{c_1}\Bigr)^{2k-2j}\frac{(c_1t-\beta)^{2k-j}(c_2t+\beta)^{j}}{j!(2k-j)!}\Biggr]$$
\end{corollary}

\paragraph{Proof.}
Formula (\ref{massimo-d}) is obtained by deriving (\ref{massimoRip-d}) and performing suitable simplifying steps, similar to those in the proof of Corollary \ref{corollarioMassimo+p}.
\hfill$\Box$
\\

If $c_1 = c_2 = c$, the sum in (\ref{massimo-d}) is canceled out and, by performing some calculation, the probability (\ref{massimo-d}) reduces to
\begin{equation}
P\{\max_{0\le s\le t}\mathcal{T}(s)\in \dif \beta\ |\ N(t) = 2k,\ V(0) = -c\} =
\end{equation}
$$  =\Biggl( \frac{(2k+1)!}{(k+1)!(k-1)!} \frac{(ct-\beta)^k(ct+ \beta)^{k-1}}{(2ct)^{2k}} +\frac{(2k+1)!}{(k+1)!k!} \frac{(c^2t^2 - \beta^2)^{k}}{(2ct)^{2k+1}}  \Biggr)\dif\beta =  $$
$$=\frac{2k+1}{2k+2} P\lbrace \max_{0\le s\le t} \mathcal{T}(s) \in \dif\beta\ |\ V(0) = -c, N(t) = 2k\rbrace\ +$$
\begin{equation}\label{massimo-dSommaPesata}
+\ \frac{1}{2k+2}P\lbrace \max_{0\le s\le t} \mathcal{T}(s) \in \dif\beta\ |\ V(0) = c, N(t) = 2k+1\rbrace
\end{equation}
which coincides with formula ($4.12$) of \cite{CO2020}. In the asymmetric motion we are not able to provide analogous relationships between the conditional distributions of the maximum.
\\

\begin{remark}
For $V(0)=-c_2,\ N(t)=3$, we have that (see formula (\ref{massimo-d}) ) for $0\le \beta \le c_1t$
\begin{equation}\label{densita3-}
P\{\max_{0\le s \le t}\mathcal{T}(t)\in \dif\beta\ |\ V(0) = -c_2,\ N(t) = 3\}/ \dif \beta\ =\ \frac{6(c_1t-\beta)(c_2t+\beta)+3\Bigl(\frac{c_2}{c_1} \Bigr)^2(c_1t-\beta)^2}{\bigl[(c_1+c_2)t\bigr]^3}
\end{equation}
$$ $$
The function (\ref{densita3-}) has a maximum at 
\begin{equation}\label{puntoDiMassimo-3}
\beta_{max} = c_1t\frac{(c_1^2-c_2^2-c_1c_2)}{2c_1-c_2} = c_1t \Bigl( 1+\frac{c_1(c_1+c_2)}{c_2^2-2c_1^2}\Bigr)
\end{equation}
The point (\ref{puntoDiMassimo-3}) is placed to the left of $\beta =c_1t$ if $c_2<\sqrt{2}c_1$
Furthermore it belongs to $(0,c_1t)$ if
\begin{equation}
c_2<\frac{\sqrt{5}-1}{2}c_1<\sqrt{2}c_1
\end{equation}
\hfill$\Diamond$
\end{remark}

In the next theorem we present the unconditional cumulative distribution function of the maximum of the initially negatively oriented telegraph process.

\begin{theorem}\label{teoremaMassimoRip-}
Let $\{\mathcal{T}(t)\}_{t\ge0}$ be the asymmetric telegraph process. For $\beta\in [0, c_1t]$
\begin{equation}\label{massimoRip-}
P\{\max_{0\le s\le t}\mathcal{T}(s)\le \beta\ |\ V(0) = -c_2\} = 
\end{equation}
$$ =e^{-\lambda t}\Biggl[\ \sum_{r=0}^\infty I_r\Bigl(\frac{2\lambda}{c_1+c_2} \sqrt{(c_1t-\beta)(c_2t+\beta)}\Bigr) \Biggl(\sqrt{\frac{c_2t+\beta}{c_1t-\beta}}\Biggr)^r+$$
$$-\sum_{r=2}^\infty I_r\Bigl(\frac{2\lambda}{c_1+c_2} \sqrt{(c_1t-\beta)(c_2t+\beta)}\Bigr)\Bigl(\frac{c_2}{c_1}\Bigr)^{r-1}\Biggl(\sqrt{\frac{c_1t-\beta}{c_2t+\beta}}\Biggr)^r\ \Biggr]$$
\end{theorem}

\paragraph{Proof.}
We begin by deriving the following joint distributions
\begin{equation}\label{massimoRip-UnitoPari}
P\{\max_{0\le s\le t}\mathcal{T}(s)\le \beta,\ \bigcup_{k=0}^\infty\{ N(t) = 2k\}\ |\ V(0) = -c_2\} = 
\end{equation}
$$ =e^{-\lambda t}\Biggl[\ \sum_{k=0}^\infty I_{2k}\Bigl(\frac{2\lambda}{c_1+c_2} \sqrt{(c_1t-\beta)(c_2t+\beta)}\Bigr) \Biggl(\sqrt{\frac{c_2t+\beta}{c_1t-\beta}}\Biggr)^{2k}+$$
$$-\sum_{k=1}^\infty I_{2k}\Bigl(\frac{2\lambda}{c_1+c_2} \sqrt{(c_1t-\beta)(c_2t+\beta)}\Bigr)\Bigl(\frac{c_2}{c_1}\Bigr)^{2k-1}\Biggl(\sqrt{\frac{c_1t-\beta}{c_2t+\beta}}\Biggr)^{2k}\ \Biggr]$$
and
\begin{equation}\label{massimoRip-UnitoDispari}
P\{\max_{0\le s\le t}\mathcal{T}(s)\le \beta,\ \bigcup_{k=0}^\infty\{ N(t) = 2k+1\}\ |\ V(0) = -c_2\} = 
\end{equation}
$$ =e^{-\lambda t}\Biggl[\ \sum_{k=0}^\infty I_{2k+1}\Bigl(\frac{2\lambda}{c_1+c_2} \sqrt{(c_1t-\beta)(c_2t+\beta)}\Bigr) \Biggl(\sqrt{\frac{c_2t+\beta}{c_1t-\beta}}\Biggr)^{2k+1}+$$
$$-\sum_{k=1}^\infty I_{2k+1}\Bigl(\frac{2\lambda}{c_1+c_2} \sqrt{(c_1t-\beta)(c_2t+\beta)}\Bigr)\Bigl(\frac{c_2}{c_1}\Bigr)^{2k}\Biggl(\sqrt{\frac{c_1t-\beta}{c_2t+\beta}}\Biggr)^{2k+1}\ \Biggr]$$
We restrict ourselves to the proof of (\ref{massimoRip-UnitoPari}) since (\ref{massimoRip-UnitoDispari}) follows in the same way.
\\In view of (\ref{massimoRip-p}) we must evaluate the following sum
$$e^{-\lambda t} \sum_{k=0}^\infty \frac{(\lambda t)^{2k}}{(2k)!}\sum_{j=0}^k  \binom{2k}{j} \frac{(c_1t-\beta)^{j}(c_2t+\beta)^{2k-j}}{\bigl[(c_1+c_2)t \bigr]^{2k}}=$$
$$ = e^{-\lambda t} \sum_{j=0}^\infty \frac{1}{j!} \sum_{k=j}^\infty \Bigl(\frac{\lambda}{c_1+c_2} \Bigr)^{2k}  \frac{(c_1t-\beta)^{j}(c_2t+\beta)^{2k-j}}{(2k-j)!} = $$
$$ = e^{-\lambda t} \sum_{j=0}^\infty \frac{1}{j!} \sum_{k=0}^\infty \Bigl(\frac{\lambda}{c_1+c_2} \Bigr)^{2k+2j}  \frac{(c_1t-\beta)^{j}(c_2t+\beta)^{2k+j}}{(2k+j)!} = $$
\begin{equation}\label{massimoRipCongiunto-PariAddendo1} =e^{-\lambda t} \sum_{k=0}^\infty I_{2k}\Bigl(\frac{2\lambda}{c_1+c_2} \sqrt{(c_1t-\beta)(c_2t+\beta)}\Bigr) \Biggl(\sqrt{\frac{c_2t+\beta}{c_1t-\beta}}\Biggr)^{2k}
\end{equation}
Again from (\ref{massimoRip-p}) we need evaluate the sum
$$-e^{-\lambda t} \sum_{k=0}^\infty \frac{(\lambda t)^{2k}}{(2k)!}\sum_{j=0}^{k-1}  \binom{2k}{j}\Bigl(\frac{c_2}{c_1}\Bigr)^{2k-1-2j} \frac{(c_1t-\beta)^{2k-j}(c_2t+\beta)^{j}}{\bigl[(c_1+c_2)t \bigr]^{2k}}=$$
$$ = -e^{-\lambda t} \sum_{j=0}^\infty \frac{1}{j!} \sum_{k=j+1}^\infty \Bigl(\frac{\lambda}{c_1+c_2} \Bigr)^{2k} \Bigl(\frac{c_2}{c_1}\Bigr)^{2k-1-2j} \frac{(c_1t-\beta)^{2k-j}(c_2t+\beta)^{j}}{(2k-j)!} = $$
$$ = -e^{-\lambda t} \sum_{j=0}^\infty \frac{1}{j!} \sum_{k=0}^\infty \Bigl(\frac{c_2}{c_1}\Bigr)^{2k+1}\Bigl(\frac{\lambda}{c_1+c_2} \Bigr)^{2k+2+2j}  \frac{(c_1t-\beta)^{2k+2+j}(c_2t+\beta)^{j}}{(2k+2+j)!} = $$
\begin{equation}\label{massimoRipCongiunto-PariAddendo2} =-e^{-\lambda t}\sum_{k=0}^\infty I_{2k+2}\Bigl(\frac{2\lambda}{c_1+c_2} \sqrt{(c_1t-\beta)(c_2t+\beta)}\Bigr)\Bigl(\frac{c_2}{c_1}\Bigr)^{2k+1}\Biggl(\sqrt{\frac{c_1t-\beta}{c_2t+\beta}}\Biggr)^{2k+2}
\end{equation}
By summing up (\ref{massimoRipCongiunto-PariAddendo1}) and (\ref{massimoRipCongiunto-PariAddendo2}) we obtain (\ref{massimoRip-UnitoPari}).
Finally, by adding formulas (\ref{massimoRip-UnitoPari}) and (\ref{massimoRip-UnitoDispari}) we arrive at the claimed result (\ref{massimoRip-}).
\hfill$\Box$
\\

The interested reader can check that 
$$ P\{\max_{0\le s\le t}\mathcal{T}(s)\le c_1t\ |\ V(0) = -c_2\}  = 1$$
\\
Furthermore
\begin{equation}\label{massimoRip-Singolarita}
P\{\max_{0\le s\le t}\mathcal{T}(s)=0\ |\ V(0) = -c_2\} = 
\end{equation}
$$= e^{-\lambda t} \Biggl[I_{0}\Bigl(\frac{2\lambda t}{c_1+c_2}\sqrt{c_1c_2}\Bigr)+\sqrt{\frac{c_2}{c_1}}I_1\Bigl(\frac{2\lambda t}{c_1+c_2}\sqrt{c_1c_2}\Bigr)\ +\ \Bigl( 1-\frac{c_1}{c_2}\Bigr)\sum_{r=2}^\infty I_r\Bigl(\frac{2\lambda t}{c_1+c_2}\sqrt{c_1c_2}\Bigr)\Biggl(\sqrt{\frac{c_2}{c_1}}\Biggr)^{r}\ \Biggr]$$
If $c_1 = \alpha c_2$, then
$$P\{\max_{0\le s\le t}\mathcal{T}(s)=0\ |\ V(0) = -c_2\} = $$
\begin{equation}\label{massimo0-Alpha}
=e^{-\lambda t} \Biggl[I_{0}\Bigl(\frac{2\lambda t \sqrt{\alpha}}{1+\alpha}\Bigr)+\frac{1}{\sqrt{\alpha}}I_{1}\Bigl(\frac{2\lambda t \sqrt{\alpha}}{1+\alpha}\Bigr) + (1-\alpha)\sum_{r=2}^\infty I_{r}\Bigl(\frac{2\lambda t \sqrt{\alpha}}{1+\alpha}\Bigr)\Bigl(\frac{1}{\sqrt{\alpha}}\Bigr)^{r}\ \Biggr]
\end{equation}
The probability mass (\ref{massimo0-Alpha}), differently from the conditional ones, depends on the product $\lambda t$. Furthermore (\ref{massimo0-Alpha}) shows that the probability of the singularity $\beta = 0$ depends on the ratio $\alpha>0$ of the two velocities only.
\\For $c_1 = c_2 = c$ (i.e. $\alpha = 1$), formula (\ref{massimoRip-Singolarita}) reduces to 
\begin{equation}\label{massimoRip-Singolarita}
P\{\max_{0\le s\le t}\mathcal{T}(s)=0\ |\ V(0) = -c\} = e^{-\lambda t} \Bigl[I_{0}\bigl(\lambda t\bigr)+I_{1}\bigl(\lambda t\bigr) \Bigr]
\end{equation}
which is independent from the velocity $c$ and coincides with ($4.18$) of \cite{CO2020}.

\begin{remark}
The following relationship shows the impact of the initial velocity on the distribution of the maximum.
By substracting (\ref{massimoRip+}) from (\ref{massimoRip-}) we have that
\begin{equation}\label{}
P\{\max_{0\le s\le t}\mathcal{T}(s)\le \beta\ |\ V(0) = -c_2\} -P\{\max_{0\le s\le t}\mathcal{T}(s)< \beta\ |\ V(0) = c_1\}=
\end{equation}
$$ = e^{-\lambda t} \Biggl[ I_0\Bigl(\frac{2\lambda}{c_1+c_2} \sqrt{(c_1t-\beta)(c_2t+\beta)}\Bigr) +\frac{c_2}{c_1}\sqrt{\frac{c_1t-\beta}{c_2t+\beta}} I_1\Bigl(\frac{2\lambda}{c_1+c_2} \sqrt{(c_1t-\beta)(c_2t+\beta)}\Bigr) +$$
$$ +\Bigl( 1-\frac{c_1}{c_2}\Bigr)\sum_{r=2}^\infty I_r\Bigl(\frac{2\lambda}{c_1+c_2} \sqrt{(c_1t-\beta)(c_2t+\beta)}\Bigr) \Biggl(\frac{c_2}{c_1}\sqrt{\frac{c_1t-\beta}{c_2t+\beta}}\Biggr)^r\ \Biggr]$$
This proves that for $c_2\ge c_1$, for all $0\le\beta\le c_1t,\ t>0$ we can state that
\begin{equation}
P\{\max_{0\le s\le t}\mathcal{T}(s)\le \beta\ |\ V(0) = -c_2\} >P\{\max_{0\le s\le t}\mathcal{T}(s)< \beta\ |\ V(0) = c_1\}
\end{equation}
\hfill$\Diamond$
\end{remark}

\section{Connection with the Euler-Poisson-Darboux equation and non-homogeneous asymmetric telegraph process}
Symmetric telegraph processes where the reversals of motion are timed by a Poisson process of rate $\lambda = \lambda(t)$, has distribution $p (x,t)$ for $|x|\le ct,\ t\ge0$, satisfying the equation
\begin{equation}\label{equazioneTelegrafo}
\frac{\partial^2 p }{\partial t^2}+2\lambda(t)  \frac{\partial p }{\partial t}= c^2\frac{\partial^2 p }{\partial x^2}
\end{equation}
with initial conditions
\begin{equation}
p(x,0) = \delta(x)\ ,\ \ \ \ p_t(x,0) = 0
\end{equation}
For $\lambda(t) = \frac{\alpha}{t}, \ \alpha>0$, equation (\ref{equazioneTelegrafo}) coincides with the Euler-Poisson-Darboux (EPD) equation and its fundamental solution reads
\begin{equation}\label{distribuzioneTelegrafoNonOmogeneoSimmetrico}
p(x,t)= \frac{1}{B\bigl(\alpha, \frac{1}{2} \bigr)} \frac{(c^2t^2-x^2)^{\alpha-1}}{(ct)^{2\alpha-1}}\ ,\ \ \ \ \ |x| <ct
\end{equation}
see \cite{FvK1992}, \cite{GO2016}. The functions of the form, for $m >0$
\begin{equation}\label{}
u(x,t)= (c^2t^2 - x^2)^{m}
\end{equation}
are themselves solutions of EPD equations, as shown in \cite{DiGO2012}, of the form
\begin{equation}\label{equazioneEPDsimmetrica}
\frac{\partial^2 u}{\partial t^2}-c^2\frac{\partial^2 u}{\partial x^2} = \frac{2m}{t}  \frac{\partial u}{\partial t}
\end{equation}
The conditional distributions of the asymmetric telegraph process, formulas (\ref{distribuzioneTelegrafoCondizionatoDispari+-}), (\ref{distribuzioneTelegrafoCondizionatoPari+}) and (\ref{distribuzioneTelegrafoCondizionatoPari-}), involve functions of the form
\begin{equation}\label{}
g(x,t) = g_{m,n}(x,t)= (c_1t -x)^m(c_2t + x)^n
\end{equation}
with $m,n>0, \ -c_2t<x<c_1t$, satifying the generalized EPD equation
\begin{equation}\label{equazioneEPDgenerale1}
\frac{\partial^2 g }{\partial t^2}-c_1c_2\frac{\partial^2 g }{\partial x^2} + (c_1-c_2)\frac{\partial^2 g }{\partial x\partial t} = \frac{1}{t}\Bigl[ (m+n) \frac{\partial p }{\partial t}+(c_1m-c_2n)\frac{\partial g }{\partial x}\Bigr]
\end{equation}
Clearly for $m=n,\ c_1 = c_1 = c$ equation (\ref{equazioneEPDgenerale1}) reduces to (\ref{equazioneEPDsimmetrica}).
\\The conditional distributions (\ref{distribuzioneTelegrafoCondizionatoDispari+-}), (\ref{distribuzioneTelegrafoCondizionatoPari+}) and (\ref{distribuzioneTelegrafoCondizionatoPari-}) have the form 
\begin{equation}\label{NucleoConXeT}
h(x,t)= C \cdot \frac{(c_1t -x)^m(c_2t + x)^n}{t^{m+n+1}}
\end{equation}
with $C$ being the normalizing constant, and satisfy the generalized  EPD equation
\begin{equation}\label{equazioneEPDgenerale2}
\frac{\partial^2 h}{\partial t^2} -c_1c_2\frac{\partial^2 h}{\partial x^2} + (c_1-c_2)\frac{\partial^2 h}{\partial x\partial t} =-\frac{1}{t}\Biggl[(m+n+2)\frac{\partial h}{\partial t} +\Bigl[ (c_1-c_2)(m+n+1)-(c_1m-c_2n)  \Bigr]\frac{\partial h}{\partial x} \Biggr]
\end{equation}
Note that (\ref{equazioneEPDgenerale2}) for $m=n=\alpha-1,\ c_1 = c_2$ becomes the classical space-symmetric one-dimensional EPD equation.
\\The conditional density functions of the maximum of the asymmetric telegraph process, (\ref{massimo+p}), (\ref{massimo+d}), (\ref{massimo-p}) and (\ref{massimo-d}), are linear combinations of functions of the form (\ref{NucleoConXeT}) with $m+n+1 = N(t)$.
\\The cumulative distribution functions of the maximum instead involve functions of the form (with $r = 0$)
\begin{equation}\label{}
k(x,t)= \frac{(c_1t -x)^m(c_2t + x)^n}{t^{m+n+r}}
\end{equation}
$r\in \mathbb{R}$, which satisfy the time-varying coefficients equations
\begin{equation}\label{equazioneEPDgenerale3}
\frac{\partial^2 k}{\partial t^2} -c_1c_2\frac{\partial^2 k}{\partial x^2} + (c_1-c_2)\frac{\partial^2 k}{\partial x\partial t} =
\end{equation}
$$ = -\frac{1}{t}\Biggl[(m+n+2r)\frac{\partial k}{\partial t} +\Bigl[ (c_1-c_2)(m+n+r)-(c_1m-c_2n)  \Bigr]\frac{\partial k}{\partial x} \Biggr]+ \frac{k}{t^2}(m+n+r)(1-r)$$
For $r = -(m+n)$ equation (\ref{equazioneEPDgenerale3}) reduces to (\ref{equazioneEPDgenerale1}). For $r = 1$ equation (\ref{equazioneEPDgenerale3}) coincides with (\ref{equazioneEPDgenerale2}). For other values of $r$ (\ref{equazioneEPDgenerale3}) is no longer an EPD equation.
\\

A random motion on the line with rightward velocity $c_1$ and leftward velocity $-c_2$ with reversals governed by a non-homogeneous Poisson process with rate $\lambda = \lambda(t)>0,\ t>0$, is described by the probabilities, $-c_2t<x<c_1t$
$$f(x,t)\dif x = P\{\mathcal{T}(t) \in \dif x, V(t) = c_1\}\ , \ \ b(x,t)\dif x = P\{\mathcal{T}(t) \in \dif x, V(t) = -c_2\}$$
which satisfy the differential system
\begin{equation}\label{sistemaDifferenzialeFB}
\begin{cases}
\frac{\partial f}{\partial t}  = -c_1\frac{\partial f}{\partial x} -\lambda(t)(f-b) \\
\frac{\partial b}{\partial t}  = c_2\frac{\partial b}{\partial x} +\lambda(t)(f-b)
\end{cases}
\end{equation}
By considering
$$ p(x,t) = f(x,t) + b(x,t) \ \ ,\ \ \ w(x, t) = f(x,t)-b(x,t)$$
the system (\ref{sistemaDifferenzialeFB}) can be reduced to the form 
\begin{equation}\label{sistemaDifferenzialePW}
\begin{cases}
\frac{\partial p}{\partial t}  = -c_1\frac{\partial f}{\partial x} +c_2\frac{\partial b}{\partial x}= -\frac{(c_1-c_2)}{2} \frac{\partial p}{\partial x}-\frac{(c_1+c_2)}{2} \frac{\partial w}{\partial x}\\
\frac{\partial w}{\partial t}  = -(c_1+c_2)\frac{\partial p}{\partial x} -2\lambda(t)w -\frac{(c_1-c_2)}{2} \frac{\partial w}{\partial x}
\end{cases}
\end{equation}
By differentiating with respect to time $t$ the first equation of (\ref{sistemaDifferenzialePW}) and subsequently replacing the second one differentiated with respect to the space variable $x$, we obtain that
\begin{equation}\label{}
\frac{\partial^2 p}{\partial t^2}  =-\frac{(c_1-c_2)}{2} \frac{\partial^2 p}{\partial x\partial t}-\frac{(c_1+c_2)}{2} \frac{\partial^2 w}{\partial x\partial t} = 
\end{equation}
$$ = -\frac{(c_1-c_2)}{2} \frac{\partial^2 p}{\partial x\partial t}-\frac{(c_1+c_2)}{2}\Bigl[ -(c_1+c_2)\frac{\partial^2 p}{\partial x^2} -2\lambda(t)\frac{\partial w}{\partial x} -\frac{(c_1-c_2)}{2} \frac{\partial^2 w}{\partial x^2}\Bigr] = $$
$$ = -\frac{(c_1-c_2)}{2} \frac{\partial^2 p}{\partial x\partial t}+\frac{(c_1+c_2)^2}{2^2}\frac{\partial^2 p}{\partial x^2} -2\lambda(t)\Bigl[ \frac{\partial p}{\partial t} +\frac{(c_1-c_2)}{2} \frac{\partial p}{\partial x} \Bigr]+$$
$$ -\frac{(c_1-c_2)}{2}\Bigl[ \frac{\partial^2 p}{\partial x\partial t} +\frac{(c_1-c_2)}{2} \frac{\partial^2 p}{\partial x^2} \Bigr] $$
Thus, we have
\begin{equation}\label{equazioneTelegrafoGenerale}
\frac{\partial^2 p}{\partial t^2}= -(c_1-c_2) \frac{\partial^2 p}{\partial x\partial t}  + c_1c_2\frac{\partial^2 p}{\partial x^2} -2\lambda(t)\Bigl[ \frac{\partial p}{\partial t} +\frac{(c_1-c_2)}{2} \frac{\partial p}{\partial x} \Bigr]
\end{equation}
Equation (\ref{equazioneTelegrafoGenerale}), for $\lambda(t) = \frac{\alpha}{t}, \ \alpha >0$, is a generalized EPD  equation which write as
\begin{equation}\label{equazioneTelegrafoGeneraleAlpha}
\frac{\partial^2 p}{\partial t^2}- c_1c_2\frac{\partial^2 p}{\partial x^2}+(c_1-c_2) \frac{\partial^2 p}{\partial x\partial t} =   -\frac{2\alpha}{t}\Bigl[ \frac{\partial p}{\partial t} +\frac{(c_1-c_2)}{2} \frac{\partial p}{\partial x} \Bigr]
\end{equation}
The fundamental solution of (\ref{equazioneTelegrafoGeneraleAlpha}) is
\begin{equation}\label{distribuzioneTelegrafoNonOmogeneo}
p(x,t)=\frac{\Gamma(2\alpha)}{\Gamma(\alpha)^2}\frac{(c_1t -x)^{\alpha-1}(c_2t+x)^{\alpha-1}}{\bigl[ (c_1+c_2)t\bigr]^{2\alpha-1}} = P\{\mathcal{T}_\alpha(t)\in \dif x\}/\dif x
\end{equation}
where $\{\mathcal{T}_\alpha(t)\}_{t\ge0}$ is the asymmetric telegraph process with rate of reversals $\lambda(t) =  \frac{\alpha}{t}, \ \alpha >0$. For $c_1 = c_2 = c$ formula (\ref{distribuzioneTelegrafoNonOmogeneo})  reduces to (\ref{distribuzioneTelegrafoNonOmogeneoSimmetrico}).
\\We note that (\ref{distribuzioneTelegrafoNonOmogeneo}) coincides with 
$$ P\{\mathcal{T}(t)\in \dif x\ |\ V(0) = c_1, \ N(t) = 2\alpha-1\} = P\{\mathcal{T}(t)\in \dif x\ |\ V(0) = -c_2, \ N(t) = 2\alpha-1\}$$
if $\alpha$ is a natural number.





\footnotesize{

}

\end{document}